\documentclass[a4paper,12pt]{amsart}

\usepackage{amsmath,amsthm,amssymb}
\usepackage[mathscr]{eucal}
 \usepackage{cite}
\usepackage{upgreek}
\usepackage[bookmarks=false]{hyperref}
\usepackage{enumerate}
\usepackage{color}

\numberwithin{equation}{section}

\theoremstyle{plain}
\newtheorem{main}{Theorem}
\newtheorem{mcor}[main]{Corollary}

\newtheorem{theorem}{Theorem}[section]

\newtheorem{lemma}[theorem]{Lemma}

\newtheorem{corollary}[theorem]{Corollary}
\theoremstyle{definition}
\newtheorem{definition}[theorem]{Definition}
\newtheorem{example}[theorem]{Example}

\newtheorem{remark}[theorem]{Remark}

\newif\ifcommentflag
\commentflagtrue%

\begin{document}

\title[Cocycle superrigidity for translation actions of product groups]{Cocycle superrigidity for translation actions \\ of product groups}

\author[Damien Gaboriau, Adrian Ioana, and Robin Tucker-Drob]{Damien Gaboriau, Adrian Ioana, and Robin Tucker-Drob}
\thanks{D.G. was supported by the C.N.R.S. and the ANR project GAMME (ANR-14-CE25-0004). \; A.I. was supported by  NSF Career Grant DMS \#1253402 and a Sloan Foundation Fellowship. }

\address{CNRS, Unit\'{e} de math\'{e}matiques Pures et Appliqu\'{e}es, ENS-Lyon, Universit\'{e} de Lyon, France}
\email{damien.gaboriau@ens-lyon.fr}
\address{Mathematics Department; University of California, San Diego, CA 90095-1555 (United States)}
\email{aioana@ucsd.edu}
\address{Texas A\&M Department of Mathematics;
Texas A\&M University, 
College Station, TX 77843-3368 (United States)}
\email{rtuckerd@math.tamu.edu}

\begin{abstract} Let $G$ be either a profinite or a connected compact group, and $\Gamma, \Lambda$ be finitely generated dense subgroups.  Assuming that the left translation action of $\Gamma$ on $G$ is strongly ergodic, we prove that any cocycle for the left-right translation action of $\Gamma\times\Lambda$ on $G$ with values in a countable group  is ``virtually" cohomologous to a group homomorphism.
Moreover, we prove that the same holds if $G$ is a (not necessarily compact) connected simple Lie group provided that $\Lambda$ contains an infinite cyclic subgroup with compact closure. We derive several applications to OE - and W$^*$- superrigidity. In particular, we obtain the first examples of compact actions of $\mathbb F_2\times\mathbb F_2$ which are W$^*$-superrigid.
 \end{abstract}

\maketitle

\section{Introduction} 

\subsection{Background}
In the last 15 years, a lot of progress has been made in the classification of measure preserving actions of countable groups on probability spaces, up to orbit equivalence  (see the surveys \cite{Sh04, Po07, Fu09, Ga10, Va10}).
Recall that two probability measure preserving ({\it p.m.p.}) actions $\Gamma\curvearrowright (X,\mu)$ and $\Lambda\curvearrowright (Y,\nu)$   are called {\it orbit equivalent} (OE) if there is an isomorphism of the probability spaces $\theta:X\rightarrow Y$ which identifies the orbits of the actions: $\theta(\Gamma x)=\Lambda\theta(x)$, for almost every $x\in X$. 
For amenable groups, the classification problem has been completely settled since the early 1980s: all ergodic p.m.p. actions of infinite amenable groups are orbit equivalent \cite{Dy59, OW80,CFW81}.

On the other hand, the non-amenable case is
far more complicated, and its study has uncovered a beautiful rigidity theory. 
In sharp contrast to the amenable case, it was shown that within certain families of non-amenable group actions, OE implies isomorphism of the groups and
 conjugacy of the actions.  
The first such rigidity phenomenon was found by R. Zimmer, in the context of actions of higher rank semisimple Lie groups, 
 by using his influential cocycle superrigidity theorem  \cite{Zi80, Zi84}. Deducing  OE rigidity from cocycle superrigidity results has since become a paradigm in the subject. Remarkably, also by using Zimmer's cocycle superrigidity theorem, A. Furman proved that ``most" ergodic p.m.p. actions $\Gamma\curvearrowright (X,\mu)$ of  higher rank lattices, including  $SL_n(\mathbb Z)\curvearrowright\mathbb T^n$ for $n\geq 3$, are  {\it OE superrigid} \cite{Fu98}. More precisely, any free  p.m.p. action $\Lambda\curvearrowright (Y,\nu)$ which is OE to $\Gamma\curvearrowright (X,\mu)$ must be (virtually) conjugate to it. Since then, numerous striking OE superrigidity results have been discovered in \cite{MS02,Po05,Po06, Ki06, Io08,PV08, Fu09, Ki09, PS09, TD14,Io14, CK15, Dr15}.

Many of these results have been obtained by applying techniques and ideas from S. Popa's deformation/rigidity theory. In all of these cases, one proves that much more than being OE superrigid, the actions in question are cocycle superrigid \cite{Po05,Po06, Io08,PV08, Fu09, PS09, TD14,Io14, Dr15}. 
The starting point was Popa's breakthrough 
 \cite{Po05} showing that any cocycle for a Bernoulli action $\Gamma\curvearrowright (X_0,\mu_0)^{\Gamma}$ of a property (T) group $\Gamma$ with values in a countable (more generally, a $\mathcal U_{\text{fin}}$) group is cohomologous to a group homomorphism. In subsequent work, Popa was able to remove the property (T) assumption on the group $\Gamma$ by requiring instead that $\Gamma$ is a product of two groups, one infinite and one non-amenable \cite{Po06}. The second author then obtained a cocycle superrigidity theorem for ergodic profinite actions of property (T) groups $\Gamma$. This shows that any cocycle for such an action with values in a countable group is virtually (i.e. after passing to an ergodic component of a finite index subgroup of $\Gamma$) cohomologous to a group homomorphism  \cite[Theorem B]{Io08}. For a generalization of this result to compact actions of property (T) groups, see \cite{Fu09}.
 
 Motivated by the analogy with Popa's cocycle superrigidity theorems \cite{Po05,Po06}, it was asked in \cite[page 340]{Io08} whether 
the conclusion of \cite[Theorem B]{Io08} holds for product groups, such as $\Gamma=\mathbb F_2\times\mathbb F_2$.  The interest in this question was especially high at the time, since a positive answer combined with N. Ozawa and S. Popa's work \cite{OP07} would have lead to the (then) first examples of  virtually W$^*$-superrigid actions. A positive answer was obtained in \cite[Theorem C]{OP08} for certain profinite actions of product groups, but only for cocycles with values into residually finite groups, and thus the desired W$^*$-superrigidity results could not be derived.
 So, while examples of W$^*$-superrigid actions were eventually found in \cite{Pe09, PV09} using different methods,  the  original question remained open.

\subsection{Cocycle superrigidity} In this paper we settle the above question. More precisely,
 we prove that profinite (and more generally, compact) actions of product groups are cocycle superrigid, under fairly general and essentially necessary conditions (see Theorem \ref{A} and Remarks \ref{rem1} and \ref{rem2}).

Before stating our main results, let us recall some terminology.
Let $\Gamma\curvearrowright (X,\mu)$ be a nonsingular action of a countable group $\Gamma$ on a standard (finite or infinite) measure space $(X,\mu)$. Let $\Delta$ be a Polish group. A measurable map $w:\Gamma\times X\rightarrow\Delta$  is called a {\it cocycle} if it satisfies the identity $w(g_1g_2,x)=w(g_1,g_2x)w(g_2,x),$ for all $g_1,g_2\in\Gamma$ and almost every $x\in X$. If $\delta:\Gamma\rightarrow\Delta$ is a homomorphism, then the formula $w(g,x)=\delta(g)$ defines a ``constant" cocycle.
 Two cocycles $w_1, w_2:\Gamma\times X\rightarrow\Delta$ are {\it cohomologous} if we can find a measurable map $\varphi:X\rightarrow\Delta$ such that $w_1(g,x)=\varphi(gx)w_2(g,x)\varphi(x)^{-1}$, for all $g\in\Gamma$ and almost every $x\in X$. 

\begin{definition}
A nonsingular action $\Gamma\curvearrowright (X,\mu)$ on a finite measure space $(X,\mu)$ is called {\it strongly ergodic} if it does not admit a non-trivial sequence of asimptotically invariant sets \cite{Sc80, CW81}. More precisely, any sequence of measurable sets $A_n\subseteq X$ such that $\mu(gA_n\Delta A_n)\rightarrow 0$, for all $g\in\Gamma$, must satisfy that  $\mu(A_n)(1-\mu(A_n))\rightarrow 0$.  
A nonsingular action $\Gamma\curvearrowright (X,\mu)$  on an infinite measure space $(X,\mu)$ is called strongly ergodic if the above condition holds when $\mu$ is replaced with an equivalent probability measure on $X$.
\end{definition}

\begin{remark}
For p.m.p. actions $\Gamma\curvearrowright^{\sigma} (X,\mu)$, strong ergodicity is implied by, and usually deduced from, {\it spectral gap}. The latter means that there is no sequence of unit vectors $F_n\in L^2(X)\ominus\mathbb C{\bf 1}$ which are almost invariant, in the sense that $\|F_n\circ\sigma_g-F_n\|_2\rightarrow 0$, for all $g\in\Gamma$. 
\end{remark}

Given a locally compact Polish group $G$, we denote by $m_G$ a fixed left invariant Haar measure. 

\begin{main}\label{A}
Let $G$ be a compact Polish group and $\Gamma,\Lambda<G$  countable dense subgroups. Consider the left-right translation action $\Gamma\times\Lambda\curvearrowright (G,m_G)$ given by $(g,h)\cdot x=gxh^{-1}$, for $g\in\Gamma,h\in\Lambda,x\in G$. Assume that $\Gamma$ acts strongly ergodically and $\Lambda$ is finitely generated. 

Let $w:(\Gamma\times\Lambda)\times G\rightarrow\Delta$ be a cocycle in a countable group $\Delta$. 

Then the following hold:

(1) Assume that $G$ is a profinite group and $\Gamma$ is finitely generated. 
Then there exists an open subgroup $G_0<G$ such that the restriction of $w$ to $(\Gamma_0\times\Lambda_0)\times G_0$ is cohomologous to a homomorphism $\delta:\Gamma_0\times\Lambda_0\rightarrow\Delta$,  where $\Gamma_0=\Gamma\cap G_0$ and $\Lambda_0=\Lambda\cap G_0$.

(2) Assume that $G$ is connected and $\tilde G$ is a simply connected locally compact group together with a continuous onto homomorphism $p:\tilde G\rightarrow G$ such that $\ker(p)<\tilde G$ is discrete. 
Then the ``lifted" cocycle $\tilde w:(\tilde\Gamma\times\tilde\Lambda)\times\tilde G\rightarrow\Delta$ given by $\tilde w((g,h),x)=w((p(g),p(h)),p(x))$ is cohomologous to a homomorphism $\delta:\tilde\Gamma\times\tilde\Lambda\rightarrow\Delta$.
\end{main}

Before giving concrete families of actions to which Theorem \ref{A} applies, let us make a few remarks regarding its hypothesis and related results.

\begin{remark}\label{rem1} Theorem \ref{A} covers general ``separately ergodic'' compact actions of product groups, under a mild strong ergodicity assumption.
Indeed, recall that a p.m.p. action $\Gamma\curvearrowright (X,\mu)$ is called {\it compact} if the closure of $\Gamma$ in the automorphism group of $(X,\mu)$ is compact. Left-right translation actions on compact groups are clearly compact.
Conversely, it can be shown that any compact p.m.p. action  of a product group $\Gamma\times\Lambda$ whose both factors $\Gamma$ and $\Lambda$ act ergodically and freely is conjugate to a left-right translation action.  
\end{remark}

\begin{remark}\label{rem2} Strong ergodicity is a necessary condition for cocycle superrigidity with countable targets. Indeed, let $\Gamma\curvearrowright^{\sigma} (X,\mu)$ be an ergodic but not strongly ergodic p.m.p. action. 
Assuming that every finite index subgroup of $\Gamma$ acts ergodically, it was shown in \cite[Corollary 27.5]{Ke10} that $\sigma$ admits cocycles into any countable group $\Delta\not=\{1\}$ which are not cohomologous to homomorphisms.
\end{remark}

\begin{remark}\label{rem} Several results in the literature establish cocycle superrigidity for general separately ergodic actions of products of locally compact Polish groups \cite{MS04, FM07, BF13}.   In these results, one assumes that the target group $\Delta$ is a group acting isometrically on a negatively curved space \cite{MS04} or a semi-simple algebraic  group \cite{FM07, BF13}.  Theorem \ref{A} establishes cocycle superrigidity for  a large, but specific class of actions of product groups. On the other hand, it imposes no additional assumptions on the target group $\Delta$, other than being countable. The first such cocycle superrigidity results, where no assumptions are placed on the target group $\Delta$, besides being countable or more generally $\mathcal U_{\text{fin}}$,  were obtained by S. Popa \cite{Po05, Po06} in the case of Bernoulli actions of property (T) and product groups.
\end{remark}

\begin{remark}\label{rem3}
Part (1) of Theorem \ref{A} generalizes \cite[Theorem C]{OP08} which derived the same conclusion under the additional assumptions that $\Gamma$ and $\Lambda$ have property ($\tau$), the action $\Gamma\times\Lambda\curvearrowright (G,m_G)$ satisfies a certain growth condition, and $\Delta$ is residually finite.
\end{remark}

Theorem \ref{A} applies to a wide class of actions. 
Indeed, there is an extensive literature providing examples of translation actions  on compact groups which have spectral gap and hence are strongly ergodic. In the next two examples, we give an overview of these results, focusing separately on the cases when $G$ is a profinite group or a connected compact Lie group.

\begin{example}  If $\Gamma$ is a subgroup of $SL_n(\mathbb Z)$, $n\geq 2$, then its congruence subgroups are given by $\Gamma(m)=\ker(\Gamma\rightarrow SL_n(\mathbb Z/m\mathbb Z))$, for $m\geq 1$.
 If  $G$ is any profinite completion of $SL_2(\mathbb Z)$ with respect to a sequence of congruence subgroups (e.g. $G=SL_2(\mathbb Z_p)$, for a prime $p$), then as a consequence of Selberg's  $3/16$ 
  theorem  \cite{Se65}, the translation action $\Gamma\curvearrowright (G,m_G)$ has spectral gap. Recently, this result has been vastly generalized, following a 2005 breakthrough of Bourgain and Gamburd \cite{BG05}.  
 Specifically, let $\Gamma<SL_n(\mathbb Z)$ be any subgroup which is Zariski dense in $SL_n$, for some $n\geq 2$, and let $G$ be any profinite completion of $\Gamma$ with respect to  a sequence of congruence subgroups (e.g. take $G$ to be the closure of $\Gamma$ in $SL_n(\mathbb Z_p)$, for a prime $p$).  Bourgain and Varj\'{u}'s work \cite{BV10} then implies that the translation action $\Gamma\curvearrowright (G,m_G)$ has spectral gap.  For the most recent developments in this direction, see \cite{SG16ab} and the references therein.

\end{example}

\begin{example}
In the early 1980s, it was shown  that for all $n\geq 2$, there exists a countable dense subgroup $\Gamma<G:=SO(n+1)$ such that the translation action $\Gamma\curvearrowright (G,m_G)$ has spectral gap (see \cite{Ma80,Su81,Dr84}). In recent years, there have been remarkable advances in extending this result to more general translation actions on connected compact Lie groups. 
Thus, Bourgain and Gamburd proved that if $\Gamma<G:=SU(d)$ is any countable dense subgroup  whose elements consist of matrices with algebraic entries, for some $d\geq 2$, then the translation action $\Gamma\curvearrowright (G,m_G)$ has spectral gap (see \cite{BG06, BG11}). More generally, Benoist and de Sacx\'{e} recently showed that the same conclusion holds if $G$ is an arbitrary compact connected simple Lie group \cite{BdS14}.

\end{example}

Part (2) of Theorem \ref{A} applies in particular if $G$ is a connected simple compact Lie group and $\Gamma,\Lambda$ are countable dense subgroups such that $\Gamma$ acts strongly ergodically and $\Lambda$ is finitely generated. Our next theorem generalizes this result to non-compact such Lie groups $G$. Moreover, it shows that the finite generation assumption on $\Lambda$ can be removed in this case.

\begin{main}\label{B} Let $G$ be a connected simple Lie group, $\Gamma,\Lambda$ be countable dense subgroups, and consider the left-right translation action $\Gamma\times\Lambda\curvearrowright (G,m_G)$. Assume that $\Gamma$ acts strongly ergodically and $\Lambda$ contains a torsion-free element $s$ such that the closure of $\langle s\rangle$
in $G$ is compact. Let $\tilde G$ be the universal cover of $G$ with the covering homomorphism $p:\tilde G\rightarrow G$.
Let $\tilde\Gamma=p^{-1}(\Gamma)$ and $\tilde\Lambda=p^{-1}(\Lambda)$.

Let $w:(\Gamma\times\Lambda)\times G\rightarrow\Delta$ be a cocycle in a countable group $\Delta$. 

Then the ``lifted" cocycle $\tilde w:(\tilde\Gamma\times\tilde\Lambda)\times\tilde G\rightarrow\Delta$ given by $\tilde w((g,h),x)=w((p(g),p(h)),p(x))$ is cohomologous to a homomorphism $\delta:\tilde\Gamma\times\tilde\Lambda\rightarrow\Delta$.
\end{main}

As we explain next, the conditions of Theorem \ref{B} are satisfied by a large class of actions. In particular, Theorem \ref{B} applies if $G=SL_2(\mathbb R)$ and $\Gamma,\Lambda$ are any countable dense subgroups with $\Gamma\subseteq SL_2(\mathbb Q)$.
\begin{example}
Let $G$ be a connected simple Lie group. Generalizing results of \cite{BG06,BG11,BdS14} in the case $G$ is compact, it was very recently shown by Boutonnet,  Salehi-Golsefidy, and the second author that the translation action $\Gamma\curvearrowright (G,m_G)$ is strongly ergodic for any countable dense subgroup $\Gamma<G$  whose elements consist of matrices with algebraic entries \cite{BISG15}.
\end{example}

\begin{remark} Let $G$ be a connected simple Lie group and denote by $A$ the set of $s\in G$ such that $\overline{\langle s\rangle}$ is compact.
 If $G=SL_2(\mathbb R)$, then $A=\{s\in G|\;|\text{Tr}(s)|<2\}\cup\{\text{Id}\}$, hence $A$ has non-empty interior. In general, $A$ has non-empty interior if and only if the rank of $G$ is equal to the rank of its maximal compact subgroups; otherwise, $A\subseteq G$ has zero Haar measure (see Lemma \ref{compact}).
 
 For a subgroup $\Lambda<G$, consider the following condition:
$(\star)$ $\Lambda$ contains a torsion-free element $s\in A$.
 By using the above result, it follows that if $G=SL_2(\mathbb R)$ or $G=Sp_{2n}(\mathbb R)$, for $n\geq 1$, then $(\star)$ holds for any countable dense subgroup $\Lambda<G$.  On the other hand, if $G=SL_n(\mathbb R)$, for $n\geq 3$, then there exist countable dense subgroups $\Lambda<G$ for which $(\star)$ fails (see Corollary \ref{centr} and Remark \ref{rem}).
\end{remark}

\subsection{Orbit equivalence superrigidity} In the rest of the introduction we present several consequences of Theorems \ref{A} and \ref{B}. We start by pointing out that the left-right translation actions considered there are OE superrigid. To this end, we first review some terminology.

\begin{definition}
Let $\Gamma\curvearrowright (X,\mu)$ and $\Delta\curvearrowright (Y,\nu)$ be free ergodic nonsingular actions.

\begin{itemize}

\item The actions are {\it conjugate} if there exist a nonsingular isomorphism $\theta:X\rightarrow Y$ and a group isomorphism $\delta:\Gamma\rightarrow\Delta$ such that $\theta(gx)=\delta(g)\theta(x)$, for all $g\in\Gamma$ and almost every $x\in X$.

\item The actions are {\it stably orbit equivalent} 
 (SOE)  if there exist measurable non-null sets $A\subseteq X$, $B\subseteq Y$, and a nonsingular isomorphism $\theta:A\rightarrow B$ such that  $\theta(\Gamma x\cap A)=\Delta\theta(x)\cap B$, for almost every $x\in A$. %If this holds for $A=X$ and $B=Y$, then the actions are called  {\it orbit equivalent} (OE).

\item We say that $\Delta\curvearrowright Y$ is {\it induced from} an action $\Delta_0\curvearrowright Y_0$ of a subgroup $\Delta_0<\Delta$ if $Y_0\subseteq Y$ is a $\Delta_0$-invariant measurable subset such that $\nu(gY_0\cap Y_0)=0$, for all $g\in\Delta\setminus\Delta_0$.
 \end{itemize}
 \end{definition}

\begin{remark}\label{induced} There are two canonical ways of producing SOE actions:  by induction and by taking quotients. To explain this, let $\Delta\curvearrowright (Y,\nu)$ be a free ergodic nonsingular action.

\begin{itemize}
\item If $\Delta\curvearrowright Y$ is induced from $\Delta_0\curvearrowright Y_0$, then the identity map  id$_{Y_0}:Y_0\rightarrow Y_0\subseteq Y$ witnesses that the actions $\Delta_0\curvearrowright Y_0$ and $\Delta\curvearrowright Y$ are SOE. 
\item If $\Sigma<\Delta$ is a normal subgroup such that the action $\Sigma\curvearrowright Y$ has a measurable fundamental domain $A\subseteq Y$, then the actions $\Delta\curvearrowright Y$ and $\Delta/\Sigma\curvearrowright Y/\Sigma$ are SOE, as witnessed by the map $A\ni y\mapsto \Sigma y\in Y/\Sigma$.
\end{itemize}
\end{remark}

\begin{remark}\label{free}
Let $G$ be a locally compact Polish group. If $g,h\in G$, then $\{x\in G|gxh^{-1}=x\}$ has positive Haar measure if and only if $g$ commutes with an open subgroup $G_0<G$ and $h=kgk^{-1}$, for some $k\in G$.  Thus, if the centralizer of every open subgroup $G_0<G$ in $G$ is trivial, then the left-right translation action $\Gamma\times\Lambda\curvearrowright (G,m_G)$ is free, for any subgroups $\Gamma,\Lambda<G$.
\end{remark}

\begin{mcor}\label{C} Let $\Gamma,\Lambda<G$ as in Theorem \ref{A} (1). Assume that any open subgroup of $G$ has trivial centralizer in $G$.
Let $\Delta\curvearrowright (Y,\nu)$ be any free ergodic nonsingular action of a countable group $\Delta$. 

Then $\Delta\curvearrowright (Y,\nu)$ is SOE to the left-right translation action $\Gamma\times\Lambda\curvearrowright (G,m_G)$ if and only if  $\Delta\curvearrowright Y$ is induced from $\Delta_0\curvearrowright Y_0$ and
 $\Delta_0\curvearrowright Y_0$ is conjugate to the left-right translation action $(\Gamma\cap G_0)\times(\Lambda\cap G_0)\curvearrowright (G_0,m_{G_0})$, for some open subgroup $G_0<G$.

Moreover,  $\Delta\curvearrowright (Y,\nu)$ is OE to  $\Gamma\times\Lambda\curvearrowright (G,m_G)$  if and only if the above holds and we additionally have that $[G:G_0]=[\Delta:\Delta_0]$.
\end{mcor}

\begin{mcor}\label{D} Let $\Gamma,\Lambda<G$ as in Theorem \ref{A}  (2) or as in Theorem \ref{B}. Consider the notation $\tilde\Gamma,\tilde\Lambda<\tilde G$ and put $Z=\ker(p)$. Assume that $\Gamma$ and $\Lambda$ contain no non-trivial central elements of $G$.  Let $\Delta\curvearrowright (Y,\nu)$ be any free ergodic nonsingular action of a countable group $\Delta$. 

Then $\Delta\curvearrowright (Y,\nu)$ is SOE to the left-right translation action $\Gamma\times\Lambda\curvearrowright (G,m_G)$ if and only if  $\Delta\curvearrowright Y$ is induced from $\Delta_0\curvearrowright Y_0$ and
 $\Delta_0\curvearrowright Y_0$ is conjugate to  $(\tilde\Gamma\times\tilde\Lambda)/\Sigma\curvearrowright \tilde G/\Sigma$, for some  subgroup $\Sigma<Z\times Z$ which contains $\{(g,g)|g\in Z\}$.

Moreover, if $\tilde G$ is compact, then $\Delta\curvearrowright (Y,\nu)$ is OE to $\Gamma\times\Lambda\curvearrowright (G,m_G)$ if and only if the above holds and we additionally have that $|\Sigma|=[\Delta:\Delta_0]\;|Z|$.
\end{mcor}

\subsection{W$^*$-superrigidity} A free ergodic p.m.p. action $\Gamma\curvearrowright (X,\mu)$ is called {\it W$^*$-superrigid} if it can be ``recovered" entirely from its group measure space von Neumann algebra, $L^{\infty}(X)\rtimes\Gamma$.
In precise terms, this means that any free ergodic p.m.p. action $\Delta\curvearrowright (Y,\nu)$  which gives rise to an isomorphic von Neumann algebra, $L^{\infty}(X)\rtimes\Gamma\cong L^{\infty}(Y)\rtimes\Delta$, must be conjugate to $\Gamma\curvearrowright (X,\mu)$. 
 
 The first families of W$^*$-superrigid actions were discovered following a surge in activity in 2009-2010.
The existence of (virtually) W$^*$-superrigid actions was established by J. Peterson in \cite{Pe09}. Shortly after, S. Popa and S. Vaes obtained the first concrete families of W$^*$-superrigid actions \cite{PV09}. The second author then showed that Bernoulli actions $\Gamma\curvearrowright (X_0,\mu_0)^{\Gamma}$ of  icc property (T) groups $\Gamma$ are  W$^*$-superrigid \cite{Io10}. Subsequently, several large classes of W$^*$-superrigid actions were found (see for instance the introduction of \cite{Dr15}).

N. Ozawa and S. Popa proved that if $\Gamma\curvearrowright (X,\mu)$ is a free ergodic p.m.p. compact action of a product of non-abelian free groups $\Gamma=\mathbb F_{n_1}\times...\times\mathbb F_{n_k}$ for $k\geq 1$, then  $L^{\infty}(X)$ is the unique Cartan subalgebra of $L^{\infty}(X)\rtimes\Gamma$, up to unitary conjugacy \cite{OP07}. 
As a consequence, any free ergodic p.m.p. action which satisfies $L^{\infty}(X)\rtimes\Gamma\cong L^{\infty}(Y)\rtimes\Lambda$ 
must be orbit equivalent to $\Gamma\curvearrowright (X,\mu)$ \cite{Si55}. 
In combination with Corollaries \ref{C} and \ref{D},  this fact leads to the first examples of W$^*$-superrigid compact actions of  $\mathbb F_2\times\mathbb F_2$.

\begin{mcor}\label{E} Let $S$ be a non-empty set of primes, and denote $K=\prod_{p\in S}PSL_2(\mathbb Z_p)$. View $PSL_2(\mathbb Z)$ as a dense subgroup of $K$, via the natural diagonal embedding. Let $\Gamma<PSL_2(\mathbb Z)$ be a finitely generated non-amenable subgroup, and denote by $G$ its closure in $K$. % ($K<G$ is open, by Strong Approximation).
 Let $\Delta\curvearrowright (Y,\nu)$ be any free ergodic p.m.p.  action of a countable group $\Delta$.

Then $L^{\infty}(G)\rtimes (\Gamma\times\Gamma)\cong L^{\infty}(Y)\rtimes\Delta$ if and only if we can find an open subgroup $G_0<G$ and a finite index subgroup $\Delta_0<\Delta$ such that $[G:G_0]=[\Delta:\Delta_0]$,
 the left-right translation action $(\Gamma\cap G_0)\times(\Gamma\cap G_0)\curvearrowright (G_0,m_{G_0})$ is conjugate to $\Delta_0\curvearrowright Y_0$, and 
 $\Delta\curvearrowright Y$ is induced from $\Delta_0\curvearrowright Y_0$.

\end{mcor}

Since $PSL_2(\mathbb Z)$ is dense in $K$, by the Strong Approximation Theorem, 
Corollary \ref{E} implies that the left-right translation action $PSL_2(\mathbb Z)\times PSL_2(\mathbb Z)\curvearrowright K$ is virtually W$^*$-superrigid, for any non-empty set of primes $S$. Moreover, since $PSL_2(\mathbb Z)\cong (\mathbb Z/2\mathbb Z)*(\mathbb Z/3\mathbb Z)$, the group $\Gamma$ can be chosen to be a non-abelian free group of arbitrary rank.

\begin{mcor}\label{F}
Let $G=SU(2)$. Let $\varphi={\arccos(\frac{a}{b})}/2$, 
for integers $a,b$ with $0<|a|<b$ and $\frac{a}{b}\not=\pm\frac{1}{2}$. Let $\Gamma<G$ be the dense subgroup generated by the matrices $A=\begin{pmatrix}e^{i\varphi} & 0 \\ 0 & e^{-i\varphi} \end{pmatrix}$ and $B=\begin{pmatrix}\cos \varphi & i\sin \varphi \\ i\sin \varphi & \cos \varphi \end{pmatrix}.$  
Let $\Delta\curvearrowright (Y,\nu)$ be any free ergodic p.m.p.  action of a countable group $\Delta$.

Then $L^{\infty}(G)\rtimes (\Gamma\times\Gamma)\cong L^{\infty}(Y)\rtimes\Delta$ if and only if the left-right translation action $\Gamma\times\Gamma\curvearrowright (G,m_G)$ is conjugate to $\Delta\curvearrowright (Y,\nu)$.
Moreover, the same holds if $\Gamma$ is replaced by any of its non-amenable subgroups.

\end{mcor}

Note that as a consequence of \cite{Sw94}, $\Gamma$ is isomorphic to $\mathbb F_2$ (see the proof of Corollary \ref{F}).

%\subsection*{Acknowledgements}

\section{Cocycle rigidity for translation actions} 
The purpose of this section is twofold. 
Firstly, we introduce and discuss basic properties of the uniform metric on the space of cocycles.
Secondly, we recall criteria for untwisting cocycles for translation actions on compact and locally compact groups.  These criteria will be used in the proofs of our main results. 

\subsection{The uniform metric on the space of cocycles}
Let $\Gamma\curvearrowright (X,\mu)$ be a nonsingular action of a countable group $\Gamma$ on a probability space $(X,\mu)$. Let $\Delta$ be a countable group.
Given two cocycles $w_1,w_2:\Gamma\times X\rightarrow\Delta$, we consider the {\it uniform distance} given by $$d_{\mu}(w_1,w_2)=\sup_{g\in\Gamma}\mu(\{x\in X|w_1(g,x)\not=w_2(g,x)\}).$$ Whenever the measure $\mu$ is clear from the context, we will use the simpler notation $d(w_1,w_2)$.

\begin{lemma}\label{average}
Let $\Gamma\curvearrowright (X,\mu)$ be a p.m.p. action of a countable group $\Gamma$. Let $w_1,w_2:\Gamma\times X\rightarrow\Delta$ be cocycles into a countable group $\Delta$. Let  $\mu=\int_{Y}\nu\text{d}\rho(\nu)$ be an integral decomposition of $\mu$, where $\rho$ is a Borel probability measure on the space $Y$ of $\Gamma$-invariant Borel probability measures on $X$.

Then $\displaystyle{d_{\mu}(w_1,w_2)\leq\int_{Y}d_{\nu}(w_1,w_2)\;\text{d}\rho(\nu)\leq 4 d_{\mu}(w_1,w_2)}.$
\end{lemma}

{\it Proof.} The first inequality is obvious. Towards proving the second inequality, we denote by $c$ the counting measure of $\Delta$. 
Given any $\Gamma$-invariant Borel probability measure $\nu$ on $X$, we define a measure preserving action
 $\Gamma\curvearrowright (X\times\Delta,\nu\times c)$ by $g\cdot (x,h)=(gx,w_1(g,x)hw_2(g,x)^{-1})$, and let
$\pi_{\nu}:\Gamma\rightarrow\mathcal U(L^2(X\times\Delta,\nu\times c))$ be the associated unitary representation. For
$\eta\in L^2(X\times\Delta,\nu\times c)$ we denote by $\|\eta\|_{2,\nu}$ the corresponding $2$-norm.

Let $\xi:X\times\Delta\rightarrow\mathbb C$ be 
the characteristic function of $X\times\{1_{\Delta}\}$. Then $\|\xi\|_{2,\nu}=1$ and \begin{equation}\label{inv1}\|\pi_{\nu}(g)\xi-\xi\|_{2,\nu}^2=2\;\nu(\{x\in X|w_1(g,x)\not=w_2(g,x)\}),\;\;\text{for every $g\in\Gamma$}.\end{equation}

From this we deduce that if $\eta$ is a $\pi_{\nu}(\Gamma)$-invariant vector, then 
\begin{equation}\label{inv2}d_{\nu}(w_1,w_2)= \frac{1}{2}\sup_{g\in\Gamma}\|\pi_{\nu}(g)\xi-\xi\|_{2,\nu}^2\leq 2\|\eta-\xi\|_{2,\nu}^2. \end{equation}

Next, by \ref{inv1} we have $\|\pi_{\mu}(g)\xi-\xi\|_{2,\mu}\leq\sqrt{2d_{\mu}(w_1,w_2)}$, for all $g\in\Gamma$.
Let $\eta\in L^2(X\times\Delta,\mu\times c)$ be the unique element of minimal $\|.\|_{2,\mu}$ in the closure of the convex  hull of $\{\pi_{\mu}(g)\xi|g\in\Gamma\}$.  Then 
$\eta$ is $\pi_{\mu}(\Gamma)$-invariant and $\|\eta-\xi\|_{2,\mu}\leq\sqrt{2d_{\mu}(w_1,w_2)}$.
Since $\mu=\int_{Y}\nu\;\text{d}\rho(\nu)$, we deduce that $\eta$ is $\pi_{\nu}(\Gamma)$-invariant, for $\rho$-almost every $\nu\in Y$. In combination with \ref{inv2} we conclude that
$$\int_{Y}d_{\nu}(w_1,w_2)\;\text{d}\rho(y)\leq \int_{Y}2\|\eta-\xi\|_{2,\nu}^2\;\text{d}\rho(\nu)=2\|\eta-\xi\|_{2,\mu}^2\leq 4d_{\mu}(w_1,w_2),$$ which finishes the proof.
\hfill$\blacksquare$

\begin{corollary}\label{average2}
Let $\Gamma\curvearrowright (X,\mu)$ be a p.m.p. action of a countable group $\Gamma$. Let $(Y,\rho)$ be a standard probability space and consider the action $\Gamma\curvearrowright (X\times Y,\mu\times\rho)$ given by $g\cdot (x,y)=(gx,y)$. Let $w_1,w_2:\Gamma\times (X\times Y)\rightarrow\Delta$ be cocycles into a countable group $\Delta$. For $y\in Y$ and $i\in\{1,2\}$, define the cocycle $w_i^y:\Gamma\times X\rightarrow\Delta$ 
by $w_i^y(g,x)=w_i(gx,y)$. 
Then $\displaystyle{\int_{Y}d_{\mu}(w_1^y,w_2^y)\;\text{d}\rho(y)\leq 4d_{\mu\times\rho}(w_1,w_2).}$
\end{corollary}

%\comment{I corrected some confusion between $\nu$ and $\rho$}
%{$\nu$ should be $\rho$ in statement and proof.}

{\it Proof.} Noticing that $\mu\times\rho=\int_{Y}(\mu\times\delta_{y})\;\text{d}\rho(y)$, the conclusion follows from Lemma \ref{average}.
\hfill$\blacksquare$

%by $w_i^y(g,x)=w_i(g,x,y)$. 
%Then $\displaystyle{\int_{Y}d_{\mu}(w_1^y,w_2^y)\;\text{d}\nu(y)\leq 4d_{\mu\times\nu}(w_1,w_2).}$
%\end{corollary}
%\comment{$\nu$ should be $\rho$ in statement and proof.}
%
%{\it Proof.} Noticing that $\mu\times\nu=\int_{Y}(\mu\times\delta_{y})\;\text{d}\nu(y)$, the conclusion follows from Lemma \ref{average}.
%\hfill$\blacksquare$
 
 \subsection{Cocycle rigidity for translation actions on compact groups} 
 Let $\Gamma$ be a countable dense subgroup of a compact group $G$, and consider the left translation action $\Gamma\curvearrowright (G,m_G)$.
If $\Gamma$ has property (T) and $G=\varprojlim\Gamma/\Gamma_n$ is a profinite completion of $\Gamma$, then the second author proved that any cocycle $w:\Gamma\times G\rightarrow\Delta$ with values into a countable group $\Delta$ is cohomologous to a cocycle which factors through the map $\Gamma\times G\rightarrow\Gamma\times \Gamma/\Gamma_n$, for some $n$ (see  \cite[Theorem B]{Io08}). 
In \cite{Fu09}, A. Furman  extended this result  to cover a more general class of compact groups, including simply connected compact groups $G$ (see \cite[Theorem 5.21]{Fu09}). 
The proofs of these results rely implicitely on the following criteria for untwisting cocycles.

\begin{theorem} \label{untwist}

Let $G$ be a compact Polish group and $\Gamma<G$ be a countable dense subgroup.
 Let $\Delta$ be a countable group and $w:\Gamma\times G\rightarrow\Delta$ be a cocycle for the left translation action $\Gamma\curvearrowright (G,m_G)$.  For $t\in G$, define a cocycle $w_t:\Gamma\times G\rightarrow\Delta$ by $w_t(g,x)=w(g,xt^{-1})$.
 Assume that for some $c\in (0,\frac{1}{32})$ we have that  
$d(w_t,w)<c$, for every $t$ in a neighborhood $V$  of the identity $1_G$ of $G$.

\begin{enumerate}
\item \emph{\cite{Io08, Fu09}} Assume that  $G$ is a profinite group. 
Then we can find an open subgroup $G_0<G$ such that  the restriction of $w$ to $\Gamma_0\times G_0$ is cohomologous to a homomorphism $\delta:\Gamma_0\rightarrow\Delta$, where $\Gamma_0=\Gamma\cap G_0$.

\item \emph{\cite{Fu09}} Assume that  $G$ is connected and $\tilde G$ is a simply connected locally compact group with a continuous onto homomorphism $p:\tilde G\rightarrow G$ such that $\ker(p)$ is discrete in $\tilde G$, and denote $\tilde\Gamma=p^{-1}(\Gamma)$.
Then the ``lifted" cocycle $\tilde w:\tilde{\Gamma}\times\tilde G\rightarrow\Delta$ for the translation action $\tilde\Gamma\curvearrowright (\tilde G, m_{\tilde G})$ given by $\tilde{w}(g,x)=w(p(g),p(x))$ is cohomologous to a homomorphism $\delta:\tilde{\Gamma}\rightarrow\Delta$. 
\end{enumerate}
\end{theorem}

{\it Proof.}  {\it (1)}.  For a proof of this part, see \cite[Theorem 3.1]{Io13}.

{\it (2)}. This part follows from \cite{Fu09}, but for the reader's convenience, we outline a proof.
Firstly, the
proof of \cite[Theorem 5.21]{Fu09} implies that for every $t\in V$, we can find a Borel map $\varphi_t:G\rightarrow\Delta$ which satisfies
$m_G(\{x\in G|\varphi_t(x) = e\}) > \frac{3}{4}$ and $w_t(g,x) = \varphi_t(gx)w(g,x)\varphi_t(x)^{-1}$, for all $g\in\Gamma$ and almost every $x\in G$ (see also the second proof of \cite[Theorem 3.1]{Io13}).
Moreover, for every $t,s\in V$ such that $ts\in V$, we have that $\varphi_{ts}(x)=\varphi_t(xs^{-1})\varphi_s(x)$, for almost every $x\in G$.

Let $W$ be a symmetric neighborhood of the identity in $\tilde G$ such that $p(W^2)\subseteq V$.
For $t\in W$, we define $\tilde{\varphi}_t:\tilde G\rightarrow\Delta$ by letting $\tilde{\varphi}_t(x)=\varphi_{p(t)}(p(x))$.
Since $p(t)\in V$, we have that
\begin{equation}\label{cocycle1} \tilde w(g,xt^{-1})=\tilde{\varphi}_t(gx)\tilde w(g,x)\tilde{\varphi}_t(x)^{-1},\;\;\text{for all $g\in\tilde\Gamma$ and almost every $x\in\tilde G$.} \end{equation}

 Moreover, if $t,s\in W$ then $p(ts)\in V$ and therefore  we have that \begin{equation}\label{cocycle2}\tilde{\varphi}_{ts}(x)=\tilde{\varphi}_t(xs^{-1})\tilde{\varphi}_s(x),\;\;\text{for almost every $x\in G$}.\end{equation}

Since $\tilde G$ is  simply connected,  the second part of the proof of \cite[Theorem 5.21]{Fu09} shows that  we can find a family of measurable maps $\{\tilde{\varphi}_t:\tilde G\rightarrow\Delta\}_{t\in\tilde G}$ which extends the family $\{\tilde{\varphi}_t:\tilde G\rightarrow\Delta\}_{t\in W}$ defined above in such a way that the identity \ref{cocycle2} holds for every $t,s\in\tilde G$. 
Then note that the set of $t\in\tilde G$ for which identity \ref{cocycle1} holds is a subgroup of $\tilde G$ which contains $W$. Since $\tilde G$ is connected, we conclude that identity \ref{cocycle1} holds for every $t\in\tilde G$.

By arguing exactly as in the end of the proof of \cite[Theorem 5.21]{Fu09} it follows that we can find a measurable map $\varphi:\tilde G\rightarrow\Delta$ such that \begin{equation}\label{cocycle3}\tilde{\varphi}_t(x)=\varphi(xt^{-1})\varphi(x)^{-1},\;\;\;\text{for almost every $(x,t)\in\tilde G\times\tilde G$}. \end{equation}

For $g\in\tilde\Gamma$, define $\psi_g:\tilde G\rightarrow\Delta$ by letting $\psi_g(x)=\varphi(gx)^{-1}\tilde w(g,x)\varphi(x)$. By combining \ref{cocycle1} and \ref{cocycle2}, we derive that for all $t\in\tilde G$ we have $\psi_g(xt^{-1})=\psi_g(x)$, for almost every $x\in\tilde G$. Fubini's theorem implies that there exists $\delta(g)\in\Delta$ such that $\psi_g(x)=\delta(g)$, for almost every $g\in\tilde G$. Since $\delta:\Gamma\rightarrow\Delta$ is clearly a homomorphism, the conclusion follows.
\hfill$\blacksquare$

 \subsection{Cocycle rigidity for translation actions on locally compact groups} 
In the proof of Theorem \ref{B} we will need the following untwisting result for cocycles of translation actions on locally compact groups. This result was obtained in \cite[Theorem 3.1]{Io14} by adapting the proof of \cite[Theorem 5.21]{Fu09} to the locally compact setting.

\begin{theorem}\emph{\cite{Io14}}\label{untwist2}
Let $G$ be a simply connected locally compact Polish group, $\Gamma<G$ a countable dense subgroup, and put $\mathcal R=\mathcal R(\Gamma\curvearrowright G)$.  
Let $\Delta$ be a countable group and $w:\mathcal R\rightarrow\Delta$ be a cocycle. Assume that there exist a measurable set $A\subseteq G$ with $0<m_G(A)<+\infty$, a constant $c\in (0,\frac{1}{32})$, and a neighborhood $V$ of $1_G$ such that for every $\alpha\in [\mathcal R|A]$ and every $t\in V$ we have $$m_G(\{x\in A|w(\alpha(x)t,xt)\not=w(\alpha(x),x)\})<c\;m_G(A).$$

Then there exist a homomorphism $\delta:\Gamma\rightarrow\Delta$ and a measurable map $\varphi: G\rightarrow\Delta$ such that we have $w(gx,x)=\varphi(gx)\delta(g)\varphi(x)^{-1}$, for all $g\in\Gamma$ and almost every $x\in G$.
\end{theorem}

\section{Main technical result}

\begin{lemma}\label{maintech} Let $G$ be a locally compact Polish group and
$\Gamma,\Lambda$ be countable subgroups. Consider the left-right translation action $\Gamma\times \Lambda \curvearrowright (G,m_G)$ given by $(g,h)\cdot x=gxh^{-1}$, for $g\in\Gamma, h\in\Lambda, x\in G$. Assume that $\Gamma$ acts strongly ergodically, $\Lambda$ is finitely generated, and $G_0=\overline{\Lambda}$ is compact.

 Let $w:(\Gamma\times\Lambda) \times G \rightarrow\Delta$ be a cocycle into a countable group $\Delta$.  
For $t\in G$, let $w_t:\Lambda\times G\rightarrow\Delta$ be the cocycle given by $w_t(h,y)=w(h,ty)$. For $x\in G$ and $t\in G_0$, let $w^x, w^x_t:\Lambda\times G_0\rightarrow\Delta$ be the cocycles for the right translation action $\Lambda\curvearrowright (G_0,m_{G_0})$ given by $w^x(h,y)=w(h,xy)$, $w^x_t(h,y)=w(h,xty)$.
%$w^x_t=w(h,xty)$.

Then for almost every $x\in G$, we have that $d(w_t^x,w^x)\rightarrow 0$, as $t\rightarrow 1_{G_0}$.

Moreover, if $G=G_0$ (i.e. $G$ is compact and $\Lambda$ is dense in $G$), then $d(w_t,w)\rightarrow 0$, as $t\rightarrow 1_G$.

\end{lemma}

Before proving Lemma \ref{maintech}, let us outline its proof, whose main idea is inspired by \cite{GTD15}.

\subsection{Outline of the proof of Lemma \ref{maintech}}
 For simplicity, assume that $G=G_0$ and put $\mu=m_G$. 
Our goal is to show that $d(w_t,w)\rightarrow 0$, as $t\rightarrow 1_G$. This is equivalent to proving that the sets $B_h^t=\{x\in G| w(h,x)=w(h,tx)\}$, for $t\in G, h\in\Lambda$, satisfy 

$$(1)\;\inf_{h\in\Lambda}\mu(B_h^t)\rightarrow 1,\;\;\text{as $t\rightarrow 1_G$}.$$

To this end, let $A_h^t = \{ x\in G | w(h,x)=w(h,xt^{-1}) \}$ and  $A_g^t = \{ x\in G | w(g,x)=w(g,xt^{-1}) \}$, for $t\in G, h\in\Lambda, g\in\Gamma$. If $U\subseteq G$ is a symmetric conjugation invariant neighborhood of $1_G$, then $\displaystyle{\int_{U}\mu(B_h^t)\;\text{d}\mu(t)=\int_{U}\mu(A_h^t)\;\text{d}\mu(t)}$, for any $h\in\Lambda$.
Using this fact, $(1)$ reduces to proving $$(2)\;\inf_{h\in\Lambda}\mu(A_h^t)\rightarrow 1,\;\;\text{as $t\rightarrow 1_G$}.$$

To prove $(2)$ we use an idea from \cite{GTD15} on how to employ strong ergodicity. Specifically, one first notices that $g^{-1}A_h^t\triangle A_h^t \subseteq G\setminus (A_g^t \cap A_g^{hth^{-1}}h)$, for any  $t\in G, h\in\Lambda, g\in\Gamma$. Thus, the sets $\{A_h^t\}$ are almost $\Gamma$-invariant, as $t\rightarrow 1_G$, uniformly in $h\in\Lambda$. Since $\Gamma$ acts strongly ergodically, it then follows that 
$$(3)\sup_{h\in\Lambda}\mu(A_h^t)(1-\mu(A_h^t))\rightarrow 0,\;\;\text{as $t\rightarrow 1_G$}.$$

Since $\mu(A_h^t)\rightarrow 1$, as $t\rightarrow 1_G$, for any $h\in\Lambda$, and $\Lambda$ is finitely generated, it is immediate that (3) implies (2).

\subsection{Proof of Lemma \ref{maintech}} Let $\nu$ be a Borel probability measure on $G/G_0$ which is $G$-quasi-invariant. We fix a Borel map $r:G/G_0\rightarrow G$ such that $r(G_0)=1_{G}$ and $r(x)G_0=x$, for every $x\in G/G_0$. We identify  $G/G_0\times G_0$ with $G$ via the Borel isomorphism $(x,y)\mapsto r(x)y$, and denote by $\mu$ the Borel probability measure on $G$ obtained by pushing forward $\nu\times m_{G_0}$ through this isomorphism. Then $\mu$ has the same null sets as $m_G$ and is invariant under the right translation action of $G_0$.

For a non-null measurable set $U\subseteq G_0$, we denote by $m_{U}$ the normalized restriction of $m_{G_0}$ to $U$.
The proof relies on the following claim.

{\bf Claim.} For every $\varepsilon>0$, there is a neighborhood $U$ of $1_{G_0}$ such that $$\int_{G/G_0\times U}d(w_t^{r(x)},w^{r(x)})\;\text{d}\nu(x)\;\text{d}m_U(t)<\varepsilon.$$ 

Before proving the claim, let us show that it implies the conclusion of the theorem.
% If $n\geq 1$, then by the claim we can find a neighborhood of $U_n$ of $1_{G_0}$ such that  $$\int_{G/G_0}\Big(\int_{U_n}d(w_t^{r(x)},w^{r(x)})\;\text{d}m_{U_n}(t)\Big)\;\text{d}\nu(x)\leq\frac{1}{2^n},\;\;\text{for all $n\geq 1$.}$$
For any $n\geq 1$, by the claim, we can find a neighborhood $U_n$ of $1_{G_0}$ such that  $$\int_{G/G_0}\Big(\int_{U_n}d(w_t^{r(x)},w^{r(x)})\;\text{d}m_{U_n}(t)\Big)\;\text{d}\nu(x)\leq\frac{1}{2^n}.$$
%\;\;\text{for all $n\geq 1$.}

This implies that the set $X$ of all $x\in G/G_0$ such that $\int_{U_n}d(w_t^{r(x)},w^{r(x)})\;\text{d}m_{U_n}(t)\rightarrow 0$, as $n\rightarrow\infty$, is co-null in $G/G_0$. 
Fix $x\in X$. If $\varepsilon>0$, then the set $B$ of $t\in G_0$ such that $d(w_t^{r(x)},w^{r(x)})<\varepsilon/2$ is non-null. Otherwise, we would have that $\int_{U_n}d(w_t^{r(x)},w^{r(x)})\;\text{d}m_{U_n}(t)\geq\varepsilon/2$, for all $n\geq 1$.
 Next, for $t_1,t_2\in B$, we have $d(w_{t_1t_2^{-1}}^{r(x)},w^{r(x)})=d(w_{t_1}^{r(x)},w_{t_2}^{r(x)})<\varepsilon$. Since $BB^{-1}$ contains a neighborhood of $1_{G_0}$ and $\varepsilon>0$ is arbitrary, we get that $d(w_t^{r(x)},w^{r(x)})\rightarrow 0$, as $t\rightarrow 1_{G_0}$. 
Now, if $z\in G_0$, then $d(w_t^{r(x)z},w^{r(x)z})=d(w_{ztz^{-1}}^{r(x)},w^{r(x)})\rightarrow 0$, as $t\rightarrow 1_{G_0}$. Since the set $\{r(x)z|x\in X, z\in G_0\}$ is co-null in $G$, we conclude that $d(w_t^{x},w^x)\rightarrow 0$, as $t\rightarrow 1_{G_0}$, for almost every $x\in G$. 
This implies the main assertion.  

For the moreover part, assume that $G_0=G$. Recall that $d(w_t^{r(x)},w^{r(x)})\rightarrow 0$, as $t\rightarrow 1_G$, for almost every $x\in G/G_0$. Since $r(G_0)=1_{G}$, $w^{1_{G}}=w$, and $w^{1_G}_t=w_t$, the moreover assertion follows.

%\comment{I don't understand why the following sentence (?)} Since $d(w_t^x,w^x)=\sup_{h\in\Lambda}m_{G_0}(\{y\in G_0 | w(h,xty)\not=w(h,xy)\}$ and $r(G_0)=1_{G}$, the conclusion of the lemma follows. 

We are therefore left with proving the claim.

{\it Proof of the claim.}
Let $S$ be a finite generating set for $\Lambda$. We may clearly assume that $\varepsilon < 1$.  Since $\Gamma$ acts strongly ergodically on $G$, we can find a finite subset $F\subseteq\Gamma$ and some $\delta >0$ such that if $A\subseteq G$ is any measurable set with $\sup _{g\in F} \mu (g^{-1}A\triangle A ) < \delta$, then either $\mu (A)< \varepsilon/4$ or $\mu (A)>1-\varepsilon/4$. Let $U\subseteq G_0$ be a symmetric, conjugation invariant, open neighborhood of $1_{G_0}$ such that every $t\in U$ satisfies
\begin{align*}
\mu (\{ x\in G | w(g,x)=w(g,xt^{-1}) \} ) &> 1-\delta /2,\;\; \text{ for all $g\in F$, and} \\
\mu (\{ x\in G | w(h,x)=w(h,xt^{-1}) \} ) &> \varepsilon/4,\;\; \text{ for all }h\in S .
\end{align*}
For $g\in\Gamma$, $h\in\Lambda$, and $t\in U$, we define the sets $A_g^t = \{ x\in G | w(g,x)=w(g,xt^{-1}) \}$ and $A_h^t = \{ x\in G | w(h,x)=w(h,xt^{-1}) \}$. Fix $t\in U$ and $h\in\Lambda$. For $g\in \Gamma$ we claim that
\begin{equation}\label{eqn:cont}
g^{-1}A_h^t\triangle A_h^t \subseteq G\setminus (A_g^t \cap A_g^{hth^{-1}}h ) .
\end{equation}
Indeed, if $x\in A_g^t \cap A_g^{hth^{-1}}h$, then $w(g,x) = w(g,xt^{-1})$ and $w(g,xh^{-1}) = w(g,xt^{-1}h^{-1} )$. Hence we have that $x\in A_h^t$ if and only if $w(h,x)=w(h,xt^{-1})$ if and only if
\begin{align*}
w(h,gx)&=w(g,xh^{-1})w(h,x)w(g,x)^{-1} = w(g,xt^{-1}h^{-1})w(h,xt^{-1})w(g,xt^{-1})^{-1} = w(h,gxt^{-1}),
\end{align*}
(the first and last equalities use the cocycle equation for $\Gamma\times \Lambda$) i.e., if and only if $gx \in A_h^t$. 

Since $t, \, hth^{-1} \in U$, for $g\in F$ we have $\mu (A_g^t)>1-\delta /2$ and $\mu (A_g^{hth^{-1}})> 1-\delta /2$. Since $\mu$ is right $G_0$-invariant, we get $\mu (A_g^t \cap A_g^{hth^{-1}}h ) > 1-\delta$. In combination with \eqref{eqn:cont}, this shows that for all $g\in F$ we have
\[
\mu (g^{-1}A_h^t\triangle A_h^t ) < \delta ,
\]
so by our choice of $\delta$ it follows that
\begin{equation}\label{eqn:either}
\text{either }\mu (A_h^t)<\varepsilon/4 \text{ or } \mu (A_h^t)>1-\varepsilon/4, \ \ \text{for all $t\in U$ and $h\in \Lambda$}  .
\end{equation}

Let $\Lambda' = \{ h\in \Lambda | \mu (A_h^t)>1-\varepsilon/4, \text{ for all }t\in U \}$. By our choice of $U$, for all $h\in S$ and $t\in U$ we have $\mu (A_h^t)>\epsilon/4$,  hence by \eqref{eqn:either} this implies that $\mu (A_h^t)>1-\varepsilon/4$. Therefore, $\Lambda'$ contains $S$. Note that $A_{h_1}^t\cap A_{h_0}^{h_1th_1^{-1}}h_1 \subseteq A_{h_0h_1}^t$, for every $h_0,h_1\in\Lambda$ and $t\in G_0$. Indeed, if $x\in A_{h_1}^t\cap A_{h_0}^{h_1th_1^{-1}}h_1$, then $w(h_1,x)=w(h_1,xt^{-1})$ and $w(h_0,xh_1^{-1}) = w(h_0, xt^{-1}h_1^{-1})$, hence
\[
w(h_0h_1, x) = w(h_0,xh_1^{-1})w(h_1,x)=w(h_0,xt^{-1}h_1^{-1})w(h_1,xt^{-1}) = w(h_0h_1,xt^{-1}),
\]
i.e., $x\in A_{h_0h_1}^t$. Therefore, if $h_0,h_1\in\Lambda'$ and $t\in U$, then $h_1th_1^{-1} \in U$ and so we deduce that $\mu (A_{h_0h_1}^t) \geq \mu (A_{h_1}^t\cap A_{h_0}^{h_1th_1^{-1}}h_1) > 1-\varepsilon/2$. Since $1-\varepsilon/2 >\varepsilon/4$, it follows from \eqref{eqn:either} that $\mu (A_{h_0h_1}^t)>1-\varepsilon/4$, whence $h_0h_1\in \Lambda'$. This shows that $\Lambda' =\Lambda$, hence we have shown that
\begin{equation}\label{uniform}
\mu (\{ x\in G|w(h,x)=w(h,xt^{-1}) \} ) >1-\varepsilon/4,\;\; \text{ for all $h\in \Lambda$ and  $t\in U$.}
\end{equation}

 Since $U$ is conjugation invariant and symmetric, by using  \ref{uniform},  for every $h\in\Lambda$ we have that 

\begin{align*}
&(\nu\times m_{G_0}\times m_U)(\{(x,y,t)\in G/G_0\times G_0\times U| w(h,r(x)y)=w(h,r(x)ty)\})=\\&(\nu\times m_{G_0}\times m_U)(\{(x,y,t)\in G/G_0\times G_0\times U|w(h,r(x)y)=w(h,r(x)yt^{-1})\})=\\&(\mu\times m_U)(\{(x,t)\in G\times U|w(h,x)=w(h,xt^{-1})\})>1-\varepsilon/4.
\end{align*}

Finally let $\Lambda$ act on $(G/G_0\times U\times G_0,\nu\times m_{U}\times m_{G_0})$ by right translations on the last component. Define the cocycles $w_1,w_2:\Lambda\times (G/G_0\times U\times G_0)\rightarrow\Delta$ by letting $w_1(h,(x,t,y))=w(h,r(x)y)$ and $w_2(h,(x,t,y))=w(h,r(x)ty)$. Then the last displayed inequality implies $d(w_1,w_2)<\varepsilon/4$. Since $w_1(h,(x,t,y))=w_t^{r(x)}(h,y)$ and $w_2(h,(x,t,y))=w^{r(x)}(h,y)$, by applying Corollary \ref{average2} we get $$\int_{G/G_0\times U}d(w_t^{r(x)},w^{r(x)})\;\text{d}\nu(x)\;\text{d}m_U(t)\leq 4d(w_1,w_2)<\varepsilon,$$ which proves the claim and the lemma.
\hfill$\blacksquare$

\section{Proof of Theorem \ref{A}} 

In this section, we prove Theorem \ref{A}.
We begin with the following elementary result.

\begin{lemma}\label{constant} Let $G$ be a locally compact Polish group and $\Gamma<G$ a countable dense subgroup. Let $\Delta$ be a countable group, $\varphi:G\rightarrow\Delta$ be a measurable map, and $\delta_1,\delta_2:\Gamma\rightarrow\Delta$ be homomorphisms such that $\varphi(gx)=\delta_1(g)\varphi(x)\delta_2(g)^{-1}$, for all $g\in\Gamma$ and almost every $x\in G$.

Then there exists an open subgroup $G_0<G$ such that $\varphi$ is constant on every left $G_0$ coset, i.e. there is a map $\tilde{\varphi}:G/G_0\rightarrow\Delta$ such that $\varphi(x)=\tilde{\varphi}(xG_0)$, for almost every $x\in G$.

In particular, if $G$ is connected, then $\varphi$ is constant, i.e. there is $h\in\Delta$ such that $\varphi(x)=h$, for almost every $x\in G$.
\end{lemma}

%\comment{Couldn't this Lemma  be stated in the more general context when is not nec. connected. The conclusion would be the existence of an open (and closed) subgroup $B$ s.t. $\varphi$ is a.s. constant on each (right-multipl) $B$-orbit ($G/B$). (same argument as in the proof of Th.A, with the invariance of $A^t=\{x\vert \varphi(x)=\varphi(xt)\}$ under the ergodic action $\Gamma$). I don't think unimodularity is relevant here. Moreover, if $G$ is compact, then $B$ can be taken to be normal; if $G$ is connected then $B=G$.\\I mean it could be explicitly stated in the more general form. Maybe at the beginning of proof of Th. A. Would make the argument common for a longer time between the various situations.}

{\it Proof.} For $t\in G$, define $A_t=\{x\in G|\varphi(xt)=\varphi(x)\}$. 
Then $A_t$ is invariant under the left translation action $\Gamma\curvearrowright (G,m_G)$. Since this action is ergodic, $A_t$ must be null or co-null, for every $t\in G$. Let $G_0<G$ be the subgroup of $t\in G$ such that $A_t$ is co-null in $G$. Since $\Delta$ is countable,  there exists a neighborhood $V$ of the identify in $G$ such that $A_t$ is non-null, for every $t\in V$. 
It follows that $V\subset G_0$ and hence $G_0$ is open in $G$.

Next, notice that the set $\{(x,t)\in G\times G_0|\varphi(x)=\varphi(xt)\}$ is co-null in $G\times G_0$. Let $B$ be the set of $x\in G$ such that $\varphi(x)=\varphi(xt)$, for almost every $t\in G_0$. Fubini's theorem implies that $B$ is co-null in $G$.
Let $x\in G$. Since $G_0$ is non-null, we can find $y\in xG_0\cap B$. Finally, if we define $\tilde\varphi(xG_0):=\varphi(y)$, then the conclusion follows.
 \hfill$\blacksquare$

\subsection{Proof of Theorem \ref{A}}
 Let $w:(\Gamma\times\Lambda)\times G\rightarrow\Delta$ be a cocycle into a countable group $\Delta$. Let $v$ be the restriction of $w$ to $\Lambda\times G$. For $t\in G$, define the cocycle $v_t:\Lambda\times G\rightarrow\Delta$ by $v_t(h,y)=v(h,ty)$. The moreover part of Lemma \ref{maintech} 
% \comment{Formally it is a Lemma ?}
  implies the existence of a neighborhood $V$ of $1_G$ such that $d(v_t,v)<\frac{1}{33}$, for any $t\in V$. We continue by proving separately the two assertions of Theorem \ref{A}.

{\it (1)} Assume that $G$ is a profinite group and $\Gamma$ is finitely generated. The first part of Theorem \ref{untwist} implies the existence of an open subgroup $G_0<G$ such that the restriction of $v$ to $\Lambda_0\times G_0$ is cohomologous to a homomorphism $\delta:\Lambda_0\rightarrow\Delta$, where $\Lambda_0=\Lambda\cap G_0$. 

Thus, there is a measurable map $\varphi:G_0\rightarrow\Delta$ such that $\varphi(xh^{-1})w(h,x)\varphi(x)^{-1}=\delta(h)$, for all $h\in\Lambda_0$ and almost every $x\in G_0$.
Let $\Gamma_0=\Gamma\cap G_0$ and define a cocycle $\tau:(\Gamma_0\times\Lambda_0)\times G_0\rightarrow\Delta$ by letting $\tau((g,h),x)=\varphi(gxh^{-1})w((g,h),x)\varphi(x)^{-1}$.

%\comment{These letters are so close, the use of $w$ and $\omega$ is confusing.}

%Let $g\in\Gamma$. 
Let $g\in\Gamma_0$. 
Since $\tau(h,x)=\delta(h)$ we deduce that for all $h\in\Lambda_0$ and almost every $x\in G_0$ we have
\begin{equation}\label{w_0} \tau(g,xh^{-1})\delta(h)=\tau((g,h),x)=\delta(h)\tau(g,x). \end{equation}

Next, since $\Gamma$ is finitely generated and $\Gamma_0<\Gamma$ is a finite index subgroup, $\Gamma_0$ admits a finite generating set $S$.
 For $g\in\Gamma_0$ and $t\in G_0$, let $A_g^t=\{x\in G_0|\tau(g,x)=\tau(g,tx)\}$.  
 By using \ref{w_0}, the beginning of the proof of Lemma \ref{constant} implies the existence of an open subgroup $G_1<G_0$ such that $A_g^t$ is co-null in $G_0$,  for every $g\in S$ and $t\in G_1$. Moreover, since $G_1<G_0$ has finite index, after replacing $G_1$ with $\cap_{g\in G_0}gG_1g^{-1}$, we may assume that $G_1$ is normal in $G_0$.

%\comment{The next argument is a bit dry!, but OK.}

%Let $\Gamma_0'$ be the set of $g\in\Gamma_0$ such that $A_g^t$ is co-null in $G_0$, for all $t\in G_1$. Then $\Gamma_0'\subseteq\Gamma_0$ is a subgroup. 
Let $\Gamma_0'$ be the set of $g\in\Gamma_0$ such that $A_g^t$ is co-null in $G_0$, for all $t\in G_1$.
%. Then $\Gamma_0'\subseteq\Gamma_0$
Then $\Gamma_0'<\Gamma_0$ is a subgroup. 
Indeed, if $g_1,g_2\in\Gamma_0'$ and $t\in G_1$, then 
we have $A_{g_1g_2}^t\supseteq g_2^{-1}A_{g_1}^{g_2tg_2^{-1}}\cap A_{g_2}^t$ 
%%%
(by the cocycle property)
%%%
and 
$A_{g_1^{-1}}^t=g_1A_{g_1}^{g_1^{-1}tg_1}$
%%%
(by the normality of $G_1$).
%%%

Since $\Gamma_0'$ contains $S$, and $S$ generates $\Gamma_0$, we conclude that $\Gamma_0'=\Gamma_0$. 

Define $\Gamma_1=\Gamma_0\cap G_1$ and $\Lambda_1=\Lambda_0\cap G_1$. Let $g\in\Gamma_1$. Since $A_g^t$ is co-null in $G_0$, for all $t\in G_1$, we get that $\tau(g,x)=\tau(g,tx)$, for almost every $(x,t)\in G_1\times G_1$.  Fubini's theorem implies that we can find $\delta(g)\in\Delta$ such that $\tau(g,x)=\delta(g)$, for almost every 
%$g\in G_1$. 
$x\in G_1$.
Hence, for every $g\in\Gamma_1$ and $h\in\Lambda_1$ we have that $\tau((g,h),x)=\delta(g)\delta(h)$, for almost every $x\in G_1$.  It follows that the restriction of $w$ to $(\Gamma_1\times\Lambda_1)\times G_1$ is cohomologous to the resulting homomorphism $\delta:\Gamma_1\times\Lambda_1\rightarrow\Delta$. 

{\it (2)} Let $\tilde w:(\tilde\Gamma\times\tilde\Lambda)\times\tilde G\rightarrow\Delta$ be the ``lift" of $w$. The second part of Theorem \ref{untwist} implies that the restriction of $\tilde w$ to $\tilde\Lambda\times\tilde G$ is cohomologous to a homomorphism $\delta:\tilde\Lambda\rightarrow\Delta$. Thus, there is a measurable map $\varphi:\tilde G\rightarrow\Delta$ such that $\varphi(xh^{-1})\tilde w(h,x)\varphi(x)^{-1}=\delta(h)$, for all $h\in\tilde\Lambda$ and almost every $x\in \tilde G$.
Define the cocycle $\tau:(\tilde\Gamma\times\tilde\Lambda)\times\tilde G\rightarrow\Delta$ by letting $\tau((g,h),x)=\varphi(gxh^{-1})w((g,h),x)\varphi(x)^{-1}$.

Let $g\in\tilde\Gamma$. Since $\tau(h,x)=\delta(h)$ we deduce that \ref{w_0} holds for every $h\in\tilde\Lambda$ and almost every $x\in G$. Lemma \ref{constant} implies that we can find $\delta(g)\in\Delta$ such that $\tau(g,x)=\delta(g)$, for almost every $x\in\tilde G$. Hence, for any $g\in\tilde \Gamma$ and $h\in\tilde\Lambda$ we have $\tau((g,h),x)=\delta(g)\delta(h)$, for almost every $x\in\tilde G$.  It follows that $\tilde w$ is cohomologous to a homomorphism $\delta:\tilde\Gamma\times\tilde\Lambda\rightarrow\Delta$. \hfill$\blacksquare$

\section{Proof of Theorem \ref{B}}

\subsection{A bounded generation lemma} Before proving Theorem \ref{B}, we establish a lemma asserting that if a dense subgroup $\Lambda$ of a simple Lie group $G$ contains an element which (topologically) generates a torus, then $G$ is locally boundedly generated by tori generated by elements in $\Lambda$. 

\begin{lemma}\label{bdgen}
Let $G$ be a connected simple Lie group of dimension $d$, and $\Lambda<G$ be a dense subgroup. Assume that $\Lambda$ contains a torsion-free element $s$ such that the closure of $\langle s\rangle$ is compact.

Then we can find $g_1,...,g_d\in\Lambda$ and $k\geq 1$ such that 
 the closure $G_i$ of $\langle g_i\rangle$ is isomorphic to $\mathbb T^{k}$,  for all $1\leq i\leq d$, and
 $\prod_{i=1}^d(G_i\cap V)$ contains a neighborhood of $1_G$, for every neighborhood $V$ of $1_G$.

\end{lemma}

{\it Proof.} Let $G_0$ be the connected component of the closure of $\langle s\rangle$. Since $G_0$ is a non-trivial connected abelian compact Lie group, it is isomorphic to $\mathbb T^{k}$, for some $k\geq 1$. Moreover, if $g\in\langle s\rangle$ is such that $\langle g\rangle=\langle s\rangle\cap G_0$, then $G_0$ is the closure of $\langle g\rangle$.

Denote by $\frak g$ and $\frak g_0$ the Lie algebras of $G$ and $G_0$, respectively. Let Ad$:G\rightarrow GL(\frak g)$ be the adjoint representation of $G$.
Since $\Lambda$ is dense in $G$, the linear span of $\text{Ad}(\Lambda)(\frak g_0)=\{\text{Ad}(h)(\eta)|h\in\Lambda,\eta\in\frak g_0\}$ is 
an 
%%non-trivial
 Ad$(G)$-invariant subspace of $\frak g$
 different from $\{0\}$. 
On the other hand, since $G$ is simple, $\frak g$ admits no non-trivial proper Ad$(G)$-invariant subspace. 
Therefore, there must exist $h_1,...,h_m\in\Lambda$ such the set  $\{\text{Ad}(h_j)(\frak g_0)|1\leq j\leq m\}$ spans $\frak g$.

For $1\leq j\leq m$, let $g_j=h_jgh_j^{-1}\in\Lambda$. Then the closure $G_j$ of $\langle g_j\rangle$ is isomorphic to $\mathbb T^k$.
Moreover, since the Lie algebra $\frak g_j$ of  $G_j$ is equal to  Ad$(h_j)(\frak g_0)$, we deduce that $\frak g$ is equal to the span of $\frak g_j, 1\leq j\leq m$.
Thus, we can find a basis $b_1,...,b_d$ of $\frak g$ such that for every $1\leq i\leq d$, $b_i$ belongs to  $\frak g_{j_i}$, for some $j_i\in\{1,...,m\}$. 

Define $\pi:\mathbb R^d\rightarrow G$ by letting $\pi(x)=\exp(x_1b_1)...\exp(x_db_d)$, for every $x=(x_i)_{1\leq i\leq d}\in\mathbb R^d$. 
Since the derivative of $\pi$ at ${\bf 0}\in\mathbb R^d$ is invertible, the inverse function theorem implies that $\pi$ is a homeomorphism between a neighborhood of ${\bf 0}$ and a neighborhood of $1_G$. Let $V$ be a neighborhood of $1_G$. Then $W=\{x=(x_i)_{1\leq i\leq d}\in\mathbb R^d|\exp(x_ib_i)\in V,\text{for every $1\leq i\leq d$}\}$ is a neighborhood of {\bf 0}. Since $\pi(W)\subseteq \prod_{i=1}^d(G_{j_i}\cap V)$, it follows that conclusion holds for $g_{j_1},...,g_{j_d}\in\Lambda$.
\hfill$\blacksquare$

\subsection{Proof of Theorem \ref{B}} Let $w:(\Gamma\times\Lambda)\times G\rightarrow\Delta$ be a cocycle into a countable group $\Delta$. 
Let $\tilde G$ be the universal cover of $G$ and $p:\tilde G\rightarrow G$ the covering homomorphism. Let $\tilde\Gamma=p^{-1}(\Gamma)$, $\tilde\Lambda=p^{-1}(\Lambda)$, and $\tilde w:(\tilde\Gamma\times\tilde\Lambda)\times\tilde G\rightarrow\Delta$ be given by $\tilde w((g,h),x)=w((p(g),p(h)),p(x))$. 

Our goal is to show that $\tilde w$ is cohomologous to a homomorphism.  
 To this end, we will analyze the restriction of $w$ to subgroups of $\Lambda$ with compact closure.
More precisely, for the next two claims, we fix a  finitely generated subgroup $\Lambda_0<\Lambda$ such that $G_0=\overline{\Lambda}_0$ is compact and connected. We denote by $\tilde G_0$ its universal cover, by  $q:\tilde G_0\rightarrow G_0$ the covering homomorphism, and put $\tilde\Lambda_0=q^{-1}(\Lambda_0)$. 

{\bf Claim 1.} 
There exist measurable maps $\varphi:G\times\tilde G_0\rightarrow\Delta$ and $\tau:(\Gamma\times\tilde\Lambda_0)\times G\rightarrow\Delta$ such that $$w(g,q(h),xq(y))=\varphi(gx,yh^{-1})\tau((g,h),x)\varphi(x,y)^{-1},$$ for all $(g,h)\in\Gamma\times\tilde\Lambda_0$ and almost every $(x,y)\in G\times\tilde G_0$.
%\comment{Again $w$ and $\omega$ is confusing}

{\it Proof of Claim 1.} For $x\in G$ and $t\in G_0$, we define cocycles $w^x, w^x_t:\Lambda_0\times G_0\rightarrow\Delta$ by letting $w^x(h,y)=w(h,xy)$ and $w^x_t(h,y)=w(h,xty)$. Since $\Gamma$ acts strongly ergodically on $G$ and $\Lambda_0$ is finitely generated, Lemma \ref{maintech} implies that $d(w_t^x,w^x)\rightarrow 0$, as $t\rightarrow 1_{G_0}$, for almost every $x\in G$. 

By applying the second part of Theorem 2.3,  we deduce that the lift of $w^x$ to $\tilde\Lambda_0\times\tilde G_0$ is cohomologous to a homomorphism. 
Thus, for almost every $x\in G$, we can find a measurable map $\varphi_x:\tilde G_0\rightarrow\Delta$ and a homomorphism $\delta_x:\tilde\Lambda_0\rightarrow\Delta$ such that \begin{equation}\label{fiber1}w(q(h),xq(y))=\varphi_x(yh^{-1})\delta_x(h)\varphi_x(y)^{-1},\;\;\text{for all $h\in\tilde\Lambda_0$ and almost every $y\in\tilde G_0$.} \end{equation}

Moreover,  we may assume that the map $G\times\tilde G_0\ni (x,y)\mapsto\varphi(x,y):=\varphi_x(y)\in\Delta$ is measurable. 
Indeed, consider the action $\tilde\Lambda_0\curvearrowright (G\times\tilde G_0,m_G\times m_{\tilde G_0})$ by right translations on the $\tilde G_0$-component, and a cocycle for this action given by $\sigma(h,(x,y))=w(q(h),xq(y))$. 
Now, the ergodic components of this action $\tilde\Lambda_0\curvearrowright G\times\tilde G_0$ are of the form $\{x\}\times\tilde G_0$, for some $x\in G$. By equation \ref{fiber1} the restriction of $\sigma$ to every such ergodic component is cohomologous to a homomorphism. By applying \cite[Theorem 3.4 and Corollary 3.11]{FMW04} it follows that $\sigma$ is cohomologous to a cocycle which does dot depend on the $\tilde G_0$-component. In other words, there is a measurable map $\varphi$ which satisfies \ref{fiber1}.

Next, we consider the action $\Gamma\times\tilde\Lambda_0\curvearrowright (G\times\tilde G_0,m_G\times m_{\tilde G_0})$ given by $(g,h)\cdot (x,y)=(gx,yh^{-1})$. Then the formula $\rho((g,h),(x,y))=w(g,q(h),xq(y))$ defines a cocycle for this action. We define $\tau:(\Gamma\times\tilde\Lambda_0)\times (G\times\tilde G_0)\rightarrow\Delta$ by letting $\tau((g,h),(x,y))=\varphi(gx, yh^{-1})^{-1}\rho((g,h),(x,y))\varphi(x,y)$. Then $\tau$ is a cocycle for the same action as above, and \ref{fiber1} gives that \begin{equation}\label{fiber2}
\tau(h,(x,y))=\delta_x(h),\;\;\text{for all $h\in\tilde\Lambda_0$ and almost every $(x,y)\in G\times\tilde G_0$.} \end{equation}

Let $g\in\Gamma$. By combining \ref{fiber2} with the cocycle equation for $\Gamma\times\tilde\Lambda_0$ we derive that for almost every $x\in G$ we have $$\tau(g,(x,yh^{-1}))\delta_x(h)=\delta_{gx}(h)\tau(g,(x,y)),\;\;\text{for all $h\in\tilde\Lambda_0$ and almost every $y\in\tilde G_0$.}$$

Since $\tilde G_0$ is connected, Lemma \ref{constant}  implies  that the map $\tilde G_0\ni y\mapsto\tau(g,(x,y))\in\Delta$ is constant, for almost every $x\in G$. Together with \ref{fiber2} we get that the map $G\times\tilde G_0\ni (x,y)\mapsto \tau((g,h),(x,y))\in\Delta$ only depends on the $x$-coordinate, almost everywhere, for all $(g,h)\in\Gamma\times\tilde\Lambda_0$, proving the claim. \hfill$\square$

 Next, let ${\mathcal R}=\mathcal R(\tilde\Gamma\curvearrowright\tilde G)$ and define the cocycle $v:{\mathcal R}\rightarrow\Delta$ by letting $v(gx,x)=\tilde w(g,x)$, for every $g\in\tilde\Gamma$ and $x\in\tilde G$.  
For a measurable set $A\subseteq\tilde G$ with $m_{\tilde G}(A)<\infty$, and $t\in\tilde G$, we define $$c(A,t)=\sup_{\alpha\in [\mathcal R|A]}m_{\tilde G}(\{x\in A|v(\alpha(x)t,xt)\not=v(\alpha(x),x)\}).$$

We next use Claim 1 to derive that $c(A,t)\rightarrow 0$, as $t\in p^{-1}(G_0)$ converges to $1_{\tilde G}$. By using a ``bounded generation" argument we will then further show  that $c(A,t)\rightarrow 0$, as $t\rightarrow 1_{\tilde G}$. This will allow us to apply Theorem \ref{untwist2} to $v$, and conclude
 that $\tilde w$ is cohomologous to a homomorphism.

{\bf Claim 2.}
 Let $A\subseteq\tilde G$ be a measurable set with $m_{\tilde G}(A)>0$ and $c>0$. Then there exists an open neighborhood $V$ of $1_{G}$ such that
 $c(A,t)\leq c$, for every  $t\in p^{-1}(G_0\cap V)$.

{\it Proof of Claim 2.} By letting $h=1_{\tilde\Lambda_0}$ in Claim 1, we deduce that for almost every $s\in\tilde G_0$ we have 

\begin{equation}\label{fiber}
w(g,xq(s))=\varphi(gx,s)\tau(g,x)\varphi(x,s)^{-1},\;\;\text{for all $g\in\Gamma$ and almost every $x\in G$.}
\end{equation}

Let $W$ be an open neighborhood of $1_{\tilde G_0}$ such that for almost every $(s_1,s_2)\in W\times W$ we have

\begin{equation}\label{s}
m_{\tilde G}(\{x\in A|\varphi(p(x),s_1)\not=\varphi(p(x),s_2)\})\leq\frac{c}{2}.
\end{equation}

Put $V=q(W)$. Then $V$ is an open neighborhood of $1_{G_0}$. Fix $\alpha\in [\tilde{\mathcal R}|A]$ and let $\gamma:A\rightarrow\tilde\Gamma$ such that $\alpha(x)=\gamma(x)x$, for almost every $x\in A$. Let $t_1,t_2\in p^{-1}(G_0\cap V)$ and let $s_1,s_2\in W$ such that $p(t_1)=q(s_1)$ and $p(t_2)=q(s_2)$. Assume
 that \ref{fiber} holds for $s\in\{s_1,s_2\}$ and \ref{s} holds for $(s_1,s_2)$. By using the definitions of $v$ and $\tilde w$ and \ref{fiber}, for every $i\in\{1,2\}$ and almost every $x\in A$ we get
 
 \begin{align*}v(\alpha(x)t_i,xt_i) &=v(\gamma(x)xt_i,xt_i)=\tilde w(\gamma(x),xt_i)=w(p(\gamma(x)),p(x)p(t_i))=w(p(\gamma(x)),p(x)q(s_i))\\&=\varphi(p(\gamma(x))p(x),s_i)\tau(p(\gamma(x)),p(x))\varphi(p(x),s_i)^{-1}\\&=\varphi(p(\alpha(x)),s_i)\tau(p(\gamma(x)),p(x))\varphi(p(x),s_i)^{-1}\end{align*}
 
 Since $\alpha$ preserves the restriction of $m_{\tilde G}$ to $A$, by combining the last identity with \ref{s} we derive that $m_{\tilde G}(\{x\in A|v(\alpha(x)t_1,xt_1)\not=v(\alpha(x)t_2,xt_2)\})\leq c.$ 
 Since our assumptions hold for almost every pair $(t_1,t_2)\in p^{-1}(G_0\cap V)\times p^{-1}(G_0\cap V)$, so does the last inequality. Since the function $(t_1,t_2)\mapsto m_{\tilde G}(\{x\in A|v(\alpha(x)t_1,xt_1)\not=v(\alpha(x)t_2,xt_2)\})$ is continuous, we therefore conclude that $m_{\tilde G}(\{x\in A|v(\alpha(x)t_1,xt_1)\not=v(\alpha(x)t_2,xt_2)\})\leq c$, for all $t_1,t_2\in p^{-1}(G_0\cap V)$ and every $\alpha\in [\mathcal R|A]$. 
 This proves the claim.
\hfill$\square$

We are now ready to prove the following:

{\bf Claim 3.} 
 Let $A\subseteq\tilde G$ be a compact set and $c>0$. Then there exists an open neighborhood $W$ of $1_{G}$ such that
 $c(A,t)\leq c$, for every  $t\in p^{-1}(W)$.
 
{\it Proof of Claim 3.}
Let $d$ be the dimension of $G$. By Lemma \ref{bdgen}, we can find $g_1,...,g_d\in\Lambda$ and $k\geq 1$ such that the closure $G_i$ of $\langle g_i\rangle$ is isomorphic to $\mathbb T^k$, for every $1\leq i\leq d$, and $\prod_{i=1}^d(G_i\cap U)$ contains a neighborhood of $1_G$, for every neighborhood $U$ of $1_G$.

Let $B$ be a compact neighborhood of $1_{\tilde G}$. Then  $AB^d$ is a compact neighborhood of $1_{\tilde G}$. By applying Claim 1,  we can find a neighborhood $V$ of $1_G$ such that  \begin{equation}\label{c(A,t)}c(AB^d,t)\leq \frac{c}{d},\;\;\text{for every $1\leq i\leq d$ and any $t\in p^{-1}(G_i\cap V)$}. \end{equation}

Let $W$ be the set of $t\in\tilde G$ that can be written as $t=t_1...t_d$, where $t_i\in B$ and $p(t_i)\in G_i\cap V$, for every $1\leq i\leq d$. Then the first paragraph guarantees that $W$ is a neighborhood of $1_{\tilde G}$.
Let $t\in W$ and consider its decomposition $t=t_1...t_d$ as above.
If $\alpha\in [\mathcal R|A]$, then we have $$\{x\in A|v(\alpha(x)t,xt)\not=v(\alpha(x),x\}\subseteq\bigcup_{i=0}^{d-1}\{x\in A|v(\alpha(x)t_1...t_{i+1},xt_1...t_{i+1})\not=v(\alpha(x)t_1...t_i,xt_1...t_i)\}.$$

Note that the map $As\ni y\mapsto \alpha(ys^{-1})s\in As$ belongs to $[\mathcal R|As]$, and that $c(A,s)\leq c(\tilde A,s)$, for every $s\in\tilde G$ and every compact set $\tilde A\subset\tilde G$ containing $A$. Indeed, the latter inequality holds because every $\alpha\in [\mathcal R|A]$ extends trivialy to an element of $[\mathcal R|\tilde A]$.
By combining these facts with the fact that $t_i\in B\cap p^{-1}(G_i\cap V)$, for all $1\leq i\leq d$, and inequality  \ref{c(A,t)}, we deduce that $$c(A,t)\leq\sum_{i=0}^{d-1} c(At_1...t_i,t_{i+1})\leq\sum_{i=0}^{d-1}c(AB^d,t_{i+1})\leq d\cdot\frac{c}{d}=c,$$
for every $t\in W$,
thus proving the claim. \hfill$\square$

Let $A\subseteq\tilde G$ be a non-null compact set. Claim 3 implies that $c(A,t)<\frac{1}{33}\; m_{\tilde G}(A)$, for any $t$ in a small enough neighborhood of $1_{\tilde G}$. Since $\tilde G$ is simply connected, Theorem \ref{untwist2} implies that the restriction of $\tilde w$ to $\tilde\Gamma\times\tilde G$ is cohomologous to a homomorphism $\delta:\tilde\Gamma\rightarrow\Delta$. By continuing as in the proof of the second part of Theorem \ref{A} it follows that $\tilde w$ is cohomologous to a homomorphim $\delta:\tilde\Gamma\times\tilde\Lambda\rightarrow\Delta$.
\hfill$\blacksquare$

\subsection{Cyclic subgroups with compact closure} \label{compactsub}
In view of Theorem \ref{B} it is natural to wonder for which Lie groups $G$ does every dense subgroup contain an infinite cyclic subgroup with compact closure. It turns out that the answer is positive if and only if the rank of $G$ is equal to the rank of its maximal compact subgroups. Equivalently, any maximal connected compact abelian subgroup (i.e. torus) of $G$ is maximal abelian.
The first step towards proving this is the following lemma.

\begin{lemma}\emph{\cite{Wi02}}\label{compact}
Let $G$ be a connected Lie group and $L<G$ a maximal compact subgroup. 
Let $A$ be the set of $g\in G$ such that $hgh^{-1}\in L$, for some $h\in G$.  
Let $T<L$ be a maximal torus. Denote by $C_G(T)$ the centralizer of $T$ in $G$, and by $C_G(T)^0$ its connected component.

\begin{enumerate}

\item If $T=C_G(T)^0$, then $A$ has non-empty interior.
\item If $T\not=C_G(T)^0$, then $A\subseteq G$ is null. 
\end{enumerate}
 \end{lemma}
 
 We are grateful to Alireza Salehi-Golsefidy for pointing out Lemma \ref{compact} to us. After proving Lemma \ref{compact}, we realized that it also follows from \cite[Theorem 1]{Wi02}. Nevertheless, we include a self-contained proof for completeness.
 
 \begin{remark}\label{cptcls}
 Note that $L$ exists and is unique up to conjugation with an element from $G$. Also, every compact subgroup of $G$ is contained in a maximal compact subgroup (for all of this, see e.g. \cite{Bo50}). Therefore, $A$ is precisely the set of $g\in G$ such that the closure of $\langle g\rangle$ is compact.

 \end{remark} 
 
 {\it Proof.} 
Since $G$ is connected, $L$ is connected (see e.g. \cite{Bo50}). By the maximal torus theorem, $T<L$ is unique up to conjugation and every element of $L$ belongs to a conjugate of $T$.

 Let $\frak g$ and $\frak l$ denote the Lie algebras of $G$ and $L$, respectively. Let Ad$:G\rightarrow GL(\frak g)$ denote the adjoint representation of $G$. Endow $\frak g$ with an Ad$(L)$-invariant inner product.
 Consider the map $\pi:G\times L\rightarrow G$ given by $\pi(g,l)=glg^{-1}$. Let $(g_0,l_0)\in G\times L$. Then it is easy to see that the rank of the derivative $\text{d}\pi(g_0,l_0):\frak g\oplus\frak l\rightarrow\frak g$ is equal to the dimension of $\{a+b-\text{Ad}(l_0)(b)|a\in\frak l, b\in\frak g\}\subseteq\frak g.$
 Thus, $\text{d}\pi(g_0,l_0)$ is surjective if and only if $\{b-\text{Ad}(l_0)(b)|b\in\frak g\ominus\frak l\}=\frak g\ominus\frak l$ if and only if there is no non-zero $b\in\frak g\ominus\frak l$ such that Ad$(l_0)(b)=b$. We are now ready to prove the two assertions.
 
 (1) Assume that $T=C_G(T)^0$. Let $l_0\in T$ such that $\langle l_0\rangle$ is dense in $T$. Let $b\in\frak g\ominus\frak l$ such that Ad$(l_0)(b)=b$. 
 Then $\{e^{xb}|x\in\mathbb R\}\subset C_G(\langle l_0\rangle)=C_G(T)$. Hence, $e^{xb}\in T\subseteq L$, for every $x\in\mathbb R$. This implies that $b\in\frak l$ and thus $b=0$. By the argument from above we deduce that $\text{d}\pi(g_0,l_0)$ is surjective, for every hence for some $g_0\in G$. Finally, the inverse function theorem implies that $A=\pi(G\times L)$ has non-empty interior.
 
 (2) Assume that $T\not=C_G(T)^0$. Then we can find an element $b$ of the Lie algebra of $C_G(T)$ which does not belong to the Lie algebra of $T$. Then $b\notin\frak l$. 
 Indeed, otherwise $\{e^{xb}|x\in\mathbb R\}\subset C_G(T)^0\cap L=T$.
   Let $l_0\in L$. Since $T<L$ is a maximal torus, we can find $g\in L$ such that $t:=gl_0g^{-1}\in T$. Then Ad$(t)(b)=0$, hence Ad$(l_0)(b')=0$, where $b'=\text{Ad}(g^{-1})(b)$. Since $b'\notin\frak l$, if $b''$ denotes the orthogonal projection of $b'$ onto $\frak g\ominus\frak l$, then $b''\not=0$ and Ad$(l_0)(b'')=0$.
   By the above argument, this implies that $\text{d}\pi(g_0,l_0)$ is not surjective, for every $g_0\in G$. As a consequence, for every $l_0\in L$ and  $g_0\in G$, there is a neighborhood $V\subseteq G\times L$ of $(g_0,l_0)$ such that $\pi(V)\subseteq G$ is contained in a proper submanifold and hence is null.
   It follows that $A=\pi(G\times L)$ is null. 
 \hfill$\blacksquare$

 \begin{corollary}\label{centr}
 Consider the notation from Lemma \ref{compact}. Moreover, assume that $G$ is a semisimple real algebraic group. 
\begin{enumerate}
\item If $C_G(T)^0=T$, then every countable  dense subgroup $\Lambda<G$ contains a torsion-free element $s$ such the closure of $\langle s\rangle$ is compact.
\item If $C_G(T)^0\not=T$, then there exists  a countable dense subgroup $\Lambda<G$ such that $\langle s\rangle$ is discrete, for every $s\in\Lambda$.  Moreover, $\Lambda$ can be chosen isomorphic to the free group $\mathbb F_2$.
\end{enumerate}
 \end{corollary}  

{\it Proof.} (1) By \cite[Theorem 1]{BG02}, $\Lambda$ contains a free subgroup $\Lambda_0$ which is still dense in $G$. 
Since $A\subset G$ has non-empty interior, any non-trivial element $s\in \Lambda_0\cap A$ satisfies the conclusion.

(2) Let $a, b$ denote free generators of the free group $\mathbb F_2$ and $w$ be a non-trivial word in
 $a$ and $b$. Consider the word map $f_{w}:G^2\rightarrow G$, where $f_w(g,h)$ is the element of $G$ obtained by replacing $a, b$ with $g, h$ in the word $w$. We claim that $f_w^{-1}(A)$ is null.
 
  By \cite[Theorem B]{Bo83} (see also \cite[Corollary 5]{La04}), the image $f_{w}(G^2)\subseteq G$ has non-empty interior. 
Let $B_w$ denote the set of all $(g,h)\in G^2$ for which the derivative $df_w(g,h):\frak g\oplus\frak g\rightarrow\frak g$ is surjective. Then $B_w$ is non-empty. Otherwise, by reasoning  as in the proof of Lemma \ref{compact} (2),  it would follow that $f_w(G^2)\subseteq G$ is null.
Since $B_w\subseteq G^2$ is Zariski-open and non-empty, we deduce that it is co-null. In particular, $f_w^{-1}(A)\cap B_w$ is co-null in $f_w^{-1}(A)$. 

If $(g,h)\in B_w$, then since $df_w(g,h)$ is surjective, we can find a neighborhood $V_{(g,h)}\subseteq G^2$ of $(g,h)$ such that $f_w^{-1}(C)\cap V_{(g,h)}\subseteq G^2$ is null, for any null set $C\subseteq G$. Since $B_w$ can be covered by countably many of the sets $\{V_{(g,h)}\}_{(g,h)\in B_w}$, we get that $f_w^{-1}(C)\cap B_w$ is null, for every null set $C\subseteq G$. Since $A$ is null, in combination with the previous paragraph, we derive that $f_w^{-1}(A)$ is null.

Finally, by \cite[Lemma 3]{Ku51}, we can find non-empty open sets $X, Y\subseteq G$ such that $\langle g,h\rangle<G$ is dense, for every $g\in X$ and $h\in Y$. 
Since $f_w^{-1}(A)$ is null for every $w\in\mathbb F_2\setminus\{1\}$, we can find $$(g,h)\in \big(X\times Y\big)\setminus\big (\bigcup_{w\in\mathbb F_2\setminus\{1\}}f_w^{-1}(A)\big).$$

Let $\Lambda<G$ be the group generated by $g$ and $h$. Then $\Lambda$ is dense in $G$. Moreover, every non-trivial element  $s\in\Lambda$ is of the form $s=f_w(g,h)$, for some $w\in\mathbb F_2\setminus\{1\}$. Thus, $s\notin A$, hence the closure of $\langle s\rangle$ is not compact (see Remark \ref{cptcls}) and therefore $\langle s\rangle$ must be discrete.
\hfill$\blacksquare$

\begin{remark} \label{rem}
It can be easily checked that  if $n\geq 2$, then the condition $T=C_G(T)^0$
holds for $G=SL_n(\mathbb R)$ if and only if $n=2$.
Moreover, the condition holds if $G=Sp_{2n}(\mathbb R)$, for $n\geq 1$, but fails if $G=SL_n(\mathbb C)$, for $n\geq 2$.

%the maximal compact subgroup of Sp_2n(R) is U(n) with the embedding A+iB-->(A B//-B A)
\end{remark}

\section{Proofs of Corollaries \ref{C} and \ref{D}}
This section is devoted to the proofs of Corollaries \ref{C} and \ref{D}. To this end, we recall from \cite[Lemma 6.2]{Io14} the following elementary result.

\begin{lemma}\label{conju}\emph{\cite{Io14}}
Let $\Gamma\curvearrowright (X,\mu)$ be a nonsingular action and $\Delta\curvearrowright (Y,\nu)$  be a free nonsingular action of countable groups $\Gamma$ and $\Delta$. Assume that 
there exist nonsingular maps $\theta:X\rightarrow Y$, $\rho:Y\rightarrow X$, and a group homomorphism $\delta:\Gamma\rightarrow\Delta$ such that $\rho(\Delta\theta(x))\subset\Gamma x$ and $\theta(gx)=\delta(g)\theta(x)$, for all $g\in\Gamma$ and almost every $x\in X$. 
Define $\Sigma=\ker(\delta)$, $\Delta_0=\delta(\Gamma)$, and $Y_0=\theta(X)$.

Then the action $\Sigma\curvearrowright X$ admits a measurable fundamental domain, the action $\Gamma/\Sigma\curvearrowright X/\Sigma$ is conjugate to $\Delta_0\curvearrowright Y_0$, and the action $\Delta\curvearrowright Y$ is induced from $\Delta_0\curvearrowright Y_0$.

Moreover, assume that $\Gamma\curvearrowright (X,\mu)$ and $\Delta\curvearrowright (Y,\nu)$ are ergodic and probability measure preserving. Then for any measurable set $C\subseteq X$ such that the restriction of $\theta$ to $C$ is one-to-one, we have that  $$[\Delta:\Delta_0]\;{\nu(\theta(C))}=|\Sigma|\;\mu(C).$$
\end{lemma}

{\it Proof.} The main assertion is precisely \cite[Lemma 6.2]{Io14}. 
To prove the moreover assertion, notice that $\theta_*\mu$ and $\nu(Y_0)^{-1}{\nu|Y_0}$ are equivalent, $\Delta_0$-invariant  probability measures on $Y_0$. Since the action $\Delta_0\curvearrowright Y_0$ is ergodic (as $\Delta\curvearrowright Y$ is ergodic), we get that $\theta_*\mu=\nu(Y_0)^{-1}{\nu|Y_0}$.  If $C\subseteq X$ is a measurable set on which $\theta$ is one-to-one, then $\theta^{-1}(\theta(C))$ is equal to the disjoint union $\bigsqcup_{g\in\Sigma}gC$. Hence, 

$$[\Delta:\Delta_0]\;\nu(\theta(C))=\nu(Y_0)^{-1}\nu(\theta(C))=\mu(\theta^{-1}(\theta(C)))=\mu(\bigsqcup_{g\in\Sigma}gC)=|\Sigma|\;\mu(C),$$
which concludes the proof of the lemma.
\hfill$\blacksquare$

\subsection{Proof of Corollary \ref{C}} Since the ``if" part of the main assertion is easy (see e.g. Remark \ref{induced}), we will only prove the ``iff" assertion. Assume that the actions $\Gamma\times\Lambda\curvearrowright (G,m_G)$ and $\Delta\curvearrowright (Y,\nu)$ are SOE. Let $A\subseteq G, B\subseteq Y$ be non-null measurable sets together with a nonsingular isomorphism $\theta:A\rightarrow B$ such that $\theta((\Gamma\times\Lambda)x\cap A)=\Delta\theta(x)\cap B$, for almost every $x\in A$.

Since the action $\Gamma\times\Lambda\curvearrowright (G,m_G)$ is ergodic, we may extend $\theta$ to a measurable map $\theta:G\rightarrow Y$ such that $\theta((\Gamma\times\Lambda)x)\subset\Delta\theta(x)$, for almost every $x\in G$. Since the action $\Delta\curvearrowright (Y,\nu)$ is free, the formula $\theta(gxh^{-1})=w((g,h),x)\theta(x)$ defines a cocycle $w:(\Gamma\times\Lambda)\times G\rightarrow\Delta$. 

By applying part (1) of Theorem \ref{A}, we can find an open subgroup $G_0<G$, a homomorphism $\delta:\Gamma_0\times\Lambda_0\rightarrow\Delta$ (where $\Gamma_0=\Gamma\cap G$ and $\Lambda_0=\Lambda\cap G_0$), and a measurable map $\varphi:G_0\rightarrow\Delta$ such that $w((g,h),x)=\varphi(gxh^{-1})\delta(g,h)\varphi(x)^{-1}$, for all $g\in\Gamma_0, h\in\Lambda_0$, and almost every $x\in G_0$. Define $\theta_0:G_0\rightarrow Y$ by letting $\theta_0(x)=\varphi(x)^{-1}\theta(x)$. Then 

\begin{equation}\label{con}\theta_0(gxh^{-1})=\delta(g,h)\theta_0(x),\;\;\text{for all $g\in\Gamma_0,h\in\Lambda_0$, and almost every $x\in G_0$.}
\end{equation}

Next, denote $\rho=\theta^{-1}:B\rightarrow A$. Since the action $\Delta\curvearrowright (Y,\nu)$ is ergodic, we may extend $\rho$ to a measurable map $\rho:Y\rightarrow G$ such that $\rho(\Delta y)\subset(\Gamma\times\Lambda)\rho(y)$, for almost every $y\in Y$. Let $q:G\rightarrow G_0$ be a Borel map such that $q(x)\in(\Gamma\times\Lambda)x$, for all $x\in G$. 

Define $\rho_0=q\circ\rho:Y\rightarrow G_0$. Then $\rho_0(\Delta\theta_0(x))\subset(\Gamma_0\times\Lambda_0)x$, for almost every $x\in G_0$.
Since $\theta_0$ and $\rho_0$ are nonsingular maps, by using \ref{con} and applying Lemma \ref{conju} we deduce that if we denote $\Sigma=\ker(\delta)$, $\Delta_0=\delta(\Gamma_0\times\Lambda_0)$, and $Y_0=\theta_0(G_0)$, then 

\begin{itemize}
\item $\Sigma\curvearrowright G_0$ admits a measurable fundamental domain, 
\item $(\Gamma_0\times\Lambda_0)/\Sigma\curvearrowright G_0/\Sigma$ is conjugate to $\Delta_0\curvearrowright Y_0$, and
\item $\Delta\curvearrowright Y$ is induced from $\Delta_0\curvearrowright Y_0$.
\end{itemize}

To finish the proof of the main assertion, it remains to prove that $\Sigma=\{(1,1)\}$.
Since $\Sigma$ preserves the probability measure $m_{G_0}$ on $G_0$, we deduce that $\Sigma$ is finite. Let $(g, h)\in\Sigma$. Since $\Sigma$ is a normal subgroup of $\Gamma_0\times\Lambda_0$, we get that $(aga^{-1}g^{-1},1)=(a,1)(g,h)(a,1)^{-1}(g,h)\in\Sigma$, for all $a\in\Gamma_0$. As $\Gamma_0<G_0$ is dense and $\Sigma$ is finite, hence closed in $G_0\times G_0$, it follows that $(aga^{-1}g^{-1},1)\in\Sigma$, for all $a\in G_0$. Since $\Sigma$ is finite, hence discrete in $G_0\times G_0$, and $G_0$ is profinite, we can find an open subgroup $G_1<G_0$ such that $aga^{-1}g^{-1}=1$, for all $a\in G_1$.  
Since the centralizer of $G_1$ in $G$ is trivial by assumption, we derive that $g=1$, and similarly that $h=1$. This proves that $\Sigma=\{(1,1)\}$.

Finally, assume that $A=G$ and $B=Y$, i.e. that $\theta:G\rightarrow Y$ is an orbit equivalence.
Let $C\subseteq G_0$ be a non-null measurable set on which $\theta_0$ is one-to-one. Then $\nu(\theta_0(C))=\nu(\theta(C))$. Since $\theta$ is measure preserving, we conclude that $\nu(\theta_0(C))=m_G(C)$. 
On the other hand, by applying the moreover part of Lemma \ref{conju} to $\theta_0$ we deduce that $[\Delta:\Delta_0]\nu(\theta_0(C))=m_G(G_0)^{-1}m_G(C)=[G:G_0]m_G(C)$. Hence, $[G:G_0]=[\Delta:\Delta_0]$, as claimed. \hfill$\blacksquare$

\subsection{Proof of Corollary \ref{D}} Since the ``if" part of the main assertion is easy (see e.g. \cite[Example 1.5. (1)]{Io14}), we will only prove the ``iff" assertion. Assume that the actions $\Gamma\times\Lambda\curvearrowright (G,m_G)$ and $\Delta\curvearrowright (Y,\nu)$ are SOE. Let $A\subseteq G, B\subseteq Y$ be non-null measurable sets together with a nonsingular isomorphism $\theta:A\rightarrow B$ such that $\theta((\Gamma\times\Lambda)x\cap A)=\Delta\theta(x)\cap B$, for almost every $x\in A$.

Since the action $\Gamma\times\Lambda\curvearrowright (G,m_G)$ is ergodic, we may extend $\theta$ to a measurable map $\theta:G\rightarrow Y$ such that $\theta((\Gamma\times\Lambda)x)\subset\Delta\theta(x)$, for almost every $x\in G$. Since the action $\Delta\curvearrowright (Y,\nu)$ is free, the formula $\theta(gxh^{-1})=w((g,h),x)\theta(x)$ defines a cocycle $w:(\Gamma\times\Lambda)\times G\rightarrow\Delta$. 

By applying part (2) of Theorem \ref{A} or Theorem \ref{B}, depending on the situation, we can find a homomorphism $\delta:\tilde\Gamma\times\tilde\Lambda\rightarrow\Delta$ and a measurable map $\varphi:\tilde G\rightarrow\Delta$ such that  for all $g\in\tilde\Gamma,h\in\tilde\Lambda$, and almost every $x\in\tilde G$ we have $w((p(g),p(h)),p(x))=\varphi(gxh^{-1})\delta(g,h)\varphi(x)^{-1}$. Define $\tilde\theta:\tilde G\rightarrow Y$ by letting $\tilde\theta(x)=\varphi(x)^{-1}\theta(p(x))$. Then 

\begin{equation}\label{conj}\tilde\theta(gxh^{-1})=\delta(g,h)\tilde\theta(x),\;\;\text{for all $g\in\tilde\Gamma,h\in\tilde\Lambda$, and almost every $x\in\tilde G$.}
\end{equation}

Next, denote $\rho=\theta^{-1}:B\rightarrow A$. Since the action $\Delta\curvearrowright (Y,\nu)$ is ergodic, we may extend $\rho$ to a measurable map $\rho:Y\rightarrow G$ such that $\rho(\Delta y)\subset(\Gamma\times\Lambda)\rho(y)$, for almost every $y\in Y$. Let $q:G\rightarrow\tilde G$ be a Borel map such that $p(q(x))=x$, for all $x\in G$. Then $q((\Gamma\times\Lambda)x)\subset(\tilde\Gamma\times\tilde\Lambda)q(x)$, for all $x\in G$.

Define $\tilde\rho=q\circ\rho:Y\rightarrow\tilde G$. It is then easy to see that $\tilde\rho(\Delta\tilde\theta(x))\subset (\tilde\Gamma\times\tilde\Lambda)x$, for almost every $x\in\tilde G$. Since both $\tilde\theta$ and $\tilde\rho$ are nonsingular maps, by using \ref{conj} and applying Lemma \ref{conju} we deduce that
if we denote $\Sigma=\ker(\delta)$, $\Delta_0=\delta(\Gamma\times\Lambda)$, and $Y_0=\tilde\theta(\tilde G)$, then 
\begin{itemize}
\item $\Sigma\curvearrowright\tilde G$ admits a measurable fundamental domain,
\item $(\tilde\Gamma\times\tilde\Lambda)/\Sigma\curvearrowright\tilde G/\Sigma$ is conjugate to $\Delta_0\curvearrowright Y_0$, and
\item $\Delta\curvearrowright Y$ is induced from $\Delta_0\curvearrowright Y_0$.
\end{itemize} 

Denote by $Z(\tilde G)$ the center of $\tilde G$.
Since $\tilde G$ is connected and $Z$ is a discrete normal subgroup, we get that $Z\subseteq Z(\tilde G)$.
Thus, if $g\in Z$, then $(g,g)\in\tilde\Gamma\times\tilde\Lambda$ acts trivially on $\tilde G$. Since the action $\Delta_0\curvearrowright Y_0$ is free, we get that $(g,g)\in\Sigma$. 
 
 To finish the proof of the main assertion, it remains to show that $\Sigma\subseteq Z\times Z$.
 Since $\Sigma\curvearrowright\tilde G$ admits a measurable fundamental domain,  $\Omega=\{g\in\tilde\Gamma|(g,1)\in\Sigma\}$ is a  subgroup of $\tilde\Gamma$ which is discrete in $\tilde G$. Let $(g,h)\in\Sigma$. Since $\Sigma$ is normal in $\tilde\Gamma\times\tilde\Lambda$, for every $a\in\tilde\Gamma$ we have that $aga^{-1}g^{-1}\in\Omega$. Since $\Omega<\tilde G$ is discrete and $\tilde G$ is connected, we get that $aga^{-1}g^{-1}=1$, for every $a\in\tilde\Gamma$. Since $\tilde\Gamma<\tilde G$ is dense, we conclude that $g\in Z(\tilde G)$. Similarly, we get that $h\in Z(\tilde G)$. Since $\tilde\Gamma\cap Z(\tilde G)=\tilde\Lambda\cap Z(\tilde G)=Z$, we derive that indeed $\Sigma\subseteq Z\times Z$.
 
 Finally, assume that $A=G$ and $B=Y$, i.e. that $\theta:G\rightarrow Y$ is an orbit equivalence. Let $C\subseteq\tilde G$ be a non-null measurable set on which $\varphi$ is constant and $p$ is one-to-one. Let $h\in\Delta$ such that $\varphi(x)=h$, for all $x\in C$. Since $\tilde\theta(x)=h^{-1}\theta(p(x))$, for all $x\in C$,  we get that $\nu(\tilde\theta(C))=\nu(\theta(p(C))$. Since $\theta$ is measure preserving and $p$ is one-to-one on $C$, we further get  $\nu(\tilde\theta(C))=m_G(p(C))=|Z|\;m_{\tilde G}(C).$ On the other hand, the moreover part of Lemma \ref{conju} implies that
$[\Delta:\Delta_0]\;\nu(\tilde\theta(C))=|\Sigma|\;m_{\tilde G}(C)$. By combining the last two identities we derive that $[\Delta:\Delta_0]\; |Z|=|\Sigma|$, as claimed.
\hfill$\blacksquare$

\section{Proofs of Corollaries \ref{E} and \ref{F}}

\subsection{Proof of Corollary \ref{E}} Put $M=L^{\infty}(G)\rtimes (\Gamma\times\Gamma)$. Note that $PSL_2(\mathbb Z)\cong(\mathbb Z/2\mathbb Z)*(\mathbb Z/3\mathbb Z)$ and its non-amenable subgroups have property (HH)$^{+}$ from \cite[Definitions 1]{OP08}. By applying \cite[Corollary A]{OP08} to the profinite action $\Gamma\times\Gamma\curvearrowright (G,m_G)$  we get that $L^{\infty}(G)$ is the unique Cartan subalgebra of $M$, up to unitary conjugacy. 
By \cite{Si55}, any free ergodic probability measure preserving action $\Delta\curvearrowright (Y,\nu)$ such that $L^{\infty}(Y)\rtimes\Delta\cong M$ is orbit equivalent to $\Gamma\times\Gamma\curvearrowright (G,m_G)$.
 
Next, as a consequence of  \cite[Theorem 1]{BV10}, the left translation action $\Gamma\curvearrowright (G,m_G)$ has spectral gap, hence it is strongly ergodic.
 Moreover, by the Strong Approximation Theorem,  $G$ is an open subgroup of $K={\prod_{p\in S}PSL_2(\mathbb Z_p)}$. Since the centralizer of any open subgroup of $K$ in $K$ is trivial, we deduce that the centralizer of any open subgroup of $G$ in $G$ is trivial. Altogether, we can apply Corollary \ref{C} to deduce the conclusion. \hfill$\blacksquare$

\subsection{Proof of Corollary \ref{F}} Let $H=SO(3)$. Then 
 $\Phi:G\rightarrow H$ given by
$$\Phi(\;\big(\begin{matrix} x & y\\-\bar{y} & \bar{x}\end{matrix}\big)\;)=\begin{pmatrix} \Re(x^2-y^2) & \Im(x^2+y^2) & -2\;\Re(xy)\\ -\Im(x^2-y^2) & \Re(x^2+y^2) &  2\;\Im(xy)\\ 2\;\Re(x\bar{y}) & 
2\;\Im(x\bar{y}) & |x|^2-|y|^2\end{pmatrix}\;\;\text{for all $x,y\in\mathbb C$ with $|x|^2+|y|^2=1$},$$
is a continuous onto homomorphism with kernel $\{\pm \text{Id}\}$
 (see e.g. \cite[Section 1.6.1]{Ha03}).

If we put $c=b^2-a^2$, then a simple computation shows that $$\Phi(A)=\begin{pmatrix}\frac{a}{b}&\frac{\sqrt{c}}{b}&0\\-\frac{\sqrt{c}}{b}&\frac{a}{b}&0\\0&0&1 \end{pmatrix}\;\;\;\text{and}\;\;\;\Phi(B)=\begin{pmatrix}1&0&0\\0&\frac{a}{b}&\frac{\sqrt{c}}{b}\\0&-\frac{\sqrt{c}}{b}&\frac{a}{b} \end{pmatrix}.$$

Since $\frac{a}{b}\notin\{0,\pm\frac{1}{2},\pm 1\}$, the main result of \cite{Sw94} implies that the subgroup of $H$ generated by $\Phi(A)$ and $\Phi(B)$ is isomorphic to $\mathbb F_2$. Therefore, $\Gamma$ is isomorphic to $\mathbb F_2$ as well. Moreover, $-\text{Id}\notin\Gamma$.  Indeed, otherwise there would be a non-trivial word $w$ in $A$ and $B$ such that $w=-\text{Id}$. But then, $\Phi(w)=\text{Id}$, or in other words the same word $w$ in $\Phi(A)$ and $\Phi(B)$ is trivial, which is a contradiction. This shows that $\Gamma$ does not contain any non-trivial central element of $G$. 

Next, put $M=L^{\infty}(G)\rtimes (\Gamma\times\Gamma)$. By applying \cite[Corollary 4.5]{OP07} to the compact action $\Gamma\times\Gamma\curvearrowright (G,m_G)$ we get that  $L^{\infty}(G)$ is the unique Cartan subalgebra of $M$, up to unitary conjugacy. 
Let $\Delta\curvearrowright (Y,\nu)$  be a free ergodic probability measure preserving action such that $L^{\infty}(Y)\rtimes\Delta\cong M$. By \cite{Si55}, this action must be orbit equivalent to $\Gamma\times\Gamma\curvearrowright (G,m_G)$.
Altogether, since $G$ is simply connected, we can apply Corollary \ref{D} and deduce that $\Delta\curvearrowright Y$ is conjugate to $\Gamma\times\Gamma\curvearrowright (G,m_G)$. 
\hfill$\blacksquare$

\end{document}